\documentclass[11pt,reqno]{amsart}

\usepackage{amsthm}
\usepackage{amssymb}
\usepackage{latexsym}
\usepackage{multicol}
\usepackage{verbatim,enumerate}
\usepackage{accents}
\usepackage[usenames]{color}
\usepackage{hyperref}
\usepackage{hyperref}
\usepackage{amsmath, amscd}
\usepackage{soul}
\usepackage{tikz}
\usetikzlibrary{matrix,arrows,decorations.pathmorphing}
\usetikzlibrary{positioning}

\advance\textwidth by 1.2in \advance\oddsidemargin by -.6in \advance\evensidemargin by -.6in

\theoremstyle{plain}
\newtheorem*{prop}{Proposition}
\newtheorem{thm}{Theorem}
\newtheorem*{lem}{Lemma}
\newtheorem*{cor}{Corollary}%[section]

\theoremstyle{definition}
\newtheorem*{example}{Example}
\newtheorem*{defn}{Definition}
\newtheorem*{theom}{Theorem}

\theoremstyle{remark}
\newtheorem*{rem}{Remark}

\newcommand{\lie}[1]{\mathfrak{#1}}

\newcommand\bc{\mathbb C}
\newcommand\bn{\mathbb N}
\newcommand\bz{\mathbb Z}

\newenvironment{pf}{\proof}{\endproof}
\newcounter{cnt}
 \makeatletter
\def\mydggeometry{\makeatletter\dg@YGRID=1\dg@XGRID=20\unitlength=0.003pt\makeatother}
\makeatother \theoremstyle{remark}

\numberwithin{equation}{section}
\makeatletter
\def\section{\def\@secnumfont{\mdseries}\@startsection{section}{1}%
  \z@{.7\linespacing\@plus\linespacing}{.5\linespacing}%
  {\normalfont\scshape\centering}}
\def\subsection{\def\@secnumfont{\bfseries}\@startsection{subsection}{2}%
  {\parindent}{.5\linespacing\@plus.7\linespacing}{-.5em}%
  {\normalfont\bfseries}}
\makeatother

\begin{document}

%%%%%%%%%%%%%%%%%%%%%%%%%%%%%%%%%%%%%%%%%%%%

\title{Root multiplicities for Borcherds algebras and graph coloring}

\author{G. Arunkumar}
\thanks{}
\address{The Institute of Mathematical Sciences, HBNI, Chennai, India}
\email{gakumar@imsc.res.in}

\author{Deniz Kus}
\address{Mathematisches Institut, Universit\" at Bonn, Germany}
\email{dkus@math.uni-bonn.de}
\thanks{D.K was partially funded under the Institutional Strategy of the University of Cologne within the German Excellence Initiative.
}
\author{R. Venkatesh}
\thanks{}
\address{Indian Institute of Technology Madras, India}
\email{rvenkatesh@iitm.ac.in}

\begin{abstract}
We establish a connection between root multiplicities for Borcherds--Kac--Moody algebras and graph coloring. We show that the generalized chromatic polynomial of the graph associated to a given Borcherds algebra can be used to give a closed formula for certain root multiplicities. Using this connection we give a second interpretation, namely that the root multiplicity of a given root coincides with the number of acyclic orientations with a unique sink of a certain graph (depending on the root). Finally, using the combinatorics of Lyndon words we construct a basis for the root spaces corresponding to these roots and determine the Hilbert series in the case when all simple roots are imaginary. As an application we give a Lie theoretic proof of Stanley's reciprocity theorem of chromatic polynomials.
\end{abstract}
\maketitle
\section{Introduction}
Borcherds algebras were introduced by R. Borcherds in \cite{Bor88} as a natural generalization of Kac--Moody algebras. The theory of Borcherds algebras gained much interest because of its application in Borcherds's proof of the Conway--Nortan conjectures on the representation theory of the Monster simple group. 
The structure theory of Borcherds algebras is very similar to the structure theory of Kac--Moody algebras; however the main point of difference is that one is allowed to have imaginary simple roots.  
The most important step in understanding the structure of these algebras is to study roots and root multiplicities; the imaginary roots being the most mysterious ones.

Effective closed formulas for the root multiplicities are unknown in general, 
except for the affine Kac--Moody algebras and some small rank Borcherds algebras; see for example \cite{LCFL14,KL13, PN02, SL07} and references therein. All these papers deal with some particular examples of small rank Borcherds algebras. The aim of the present paper is to give a combinatorial formula for certain root multiplicities using tools from algebraic graph theory.

Let $\lie g$ be a Borcherds algebra with simple roots $\{\alpha_i: i\in I\}$, where $I$ is a finite or countably infinite set and let $G$ be the graph of $\lie g$ with vertex set $I$. We fix a tuple of non--negative integers $\mathbf{k}=(k_i : i\in I)$ with finitely many non--zero entries. Using the celebrated denominator identity, we show that certain root multiplicities have a close relationship with the generalized chromatic polynomial of $G$. 
The chromatic polynomial is an important tool in graph theory and was introduced by Birkhoff to attack the four color conjecture. It counts the number of graph colorings as a function of the number of colors. 

 A multicoloring of $G$ associated to $\bold k$ is an assignment of colors to the vertices of $G$ in which each vertex $i\in I$ receives exactly $k_i$ colors such that adjacent vertices receive only disjoint colors. The generalized chromatic polynomial counts the number of distinct proper vertex multicolorings of $G$ and is denoted by $\pi^G_{\bold k}(q)$. The following statement is one of the main results of this paper, which gives a Lie theoretic interpretation of the generalized chromatic polynomial; see Section~\ref{section3} for more details. 

\begin{theom}
 	Let $G$ be the graph of a given Borcherds algebra $\lie g$. Further assume that $k_i\leq 1$ for all real simple roots $\alpha_i$ and let $\eta(\bold k)=\sum\limits_{i\in I}k_i\alpha_i$ be a root of $\lie g$. 
 	Then,
 	$$\pi^G_{\bold k}(q)=(-1)^{\sum\limits_{i\in I}k_i} \left(\prod\limits_{\alpha \in \Delta_+} (1 - e^{-\alpha})^{q\dim \mathfrak{g}_{\alpha}}\right)[e^{-\eta(\bold k)}]$$ 
 	where $\Delta_+$ denotes the set of positive roots of $\lie g$ and $f[e^{\eta}]$ denotes the coefficient of $e^{\eta}$ in $f$. 
	\end{theom}
 We remark that the above theorem is a generalization of \cite[Theorem 1.1]{VV15} where the authors considered the special case when $\lie g$ is a Kac--Moody algebra and $k_i=1$ for all $i\in I$.
As an immediate application of the above theorem we obtain the following combinatorial formula.
\begin{theom}
Under the assumptions of the above theorem we have,
	$$ \text{ mult } \eta(\bold k) = \sum\limits_{\ell | \bold k}\frac{\mu(\ell)}{\ell}\ |\pi^G_{\bold k/\ell}(q)[q]|$$
	 where $|\pi^G_{\bold k}(q)[q]|$ denotes the absolute value of the coefficient of $q$ in $\pi^G_{\bold k}(q)$ and $\mu$ is the M\"{o}bius function.

\end{theom}
The generalized chromatic polynomial can be easily computed for many families of graphs, 
including complete graphs, trees or more generally chordal graphs and hence the above formula is more effective in these cases; 
see Section~\ref{examples} for some explicit examples. In the case when $G$ is a complete graph the above formula reduces to the well--known formula of Witt \cite{Wi37}. For other types of generalizations of Witt's formula we refer to \cite{KK96}.

Another useful interpretation of the root multiplicities is in terms of acyclic orientations. For the rest of the discussion, we consider the subgraph of $G$ spanned by the vertices $i\in I$ such that $k_i\neq 0$ which we also denote by $G$ for simplicity. For a node $i\in I$, we denote by $O_{i}(G)$ the set of acyclic orientations of $G$ with unique sink $i$. Using a result of Greene and Zaslavsky \cite{GZ83} we prove that
$$\text{ mult } \eta(\bold k) = \sum\limits_{\ell | \bold k}\frac{\mu(\ell)}{\ell}\frac{1}{(\bold k/\ell) !}\ O_{i}(G(\bold k/\ell))$$
where $G(\bold k)$ is a graph called the join of $G$ with respect to $\bold k$; see Section~\ref{clique} for a precise definition. 

In the second part of this paper we will use the combinatorial formula and the combinatorics of Lyndon words to construct a basis for the root space $\lie g_{\eta(\bold k)}$ associated to the root $\eta(\bold k).$ Indeed we will give a family of bases; one for each fixed index $i\in I$. We construct a natural combinatorial model $C^i(\bold k,G)$ and associate to each element $\bold w\in C^i(\bold k,G)$ a Lie word $\iota(\bold w)\in \lie g_{\eta(\bold k)}$. The theorem we prove is the following; see Section~\ref{section4} for more details.
\begin{theom}
The set $\{\iota(\bold w): \bold w\in C^i(\bold k,G)\}$ is a basis of $\lie g_{\eta(\bold k)}$.
\end{theom}  

The most interesting case arises when all the simple roots of $\lie g$ are imaginary. 
In this case the positive part of $\lie g$ becomes isomorphic to the free partially commutative Lie algebra associated with $G$ which is the Lie algebra freely generated by the elements $e_i, i\in I$ and defining relations $[e_i, e_j]=0$ whenever there is no edge between the vertices $i$ and $j$.  For example, it was proved in \cite{PS06} that the free partially commutative Lie algebra of $G$ appears naturally as the associated graded Lie algebra of the right angled Artin group of $G$ which is defined using lower central series; see Section~\ref{section5} for more details.

Finding a basis for these types of Lie algebras was the central study of the papers \cite{CM13},\cite{Chi06} \cite{DK92},\cite{Ly54} and \cite{PO11}; the best known such basis is the Lyndon--Shirshov basis or Chibrikov's right--normed basis for free Lie algebras. Our methods of finding a basis for free partially commutative Lie algebras is different and also the basis we obtain. 

Although the definition of the generalized chromatic polynomial depends on $q\in \mathbb N$, 
it makes sense to consider the evaluation of that polynomial at any integer. 
In the last part of this paper, we show that these numbers show up in the Hilbert series of the $q$--fold tensor product of the universal enveloping algebra associated to free partially commutative Lie algebras. We use this interpretation to give a simple Lie theoretic proof of Stanley's reciprocity theorem of chromatic polynomials 
\cite{S73}. To explain this further, let $\mathcal{O}$ be an acyclic orientation of $G$. 
For a map $\sigma: I\to \{1, 2,\dots,q\}$, we say that $(\sigma, \mathcal{O})$ is a $q$--compatible pair if for each directed edge $i\rightarrow j$ in $\mathcal{O}$ we have $\sigma(i)\ge \sigma(j)$. 
Stanley proved that the number of $q$--compatible pairs of $G$ is (up to a sign) equal to $\pi_{\bold 1}(-q)$ where $\bold 1=(k_i=1 : i\in I)$. This result will be an immediate consequence from the following theorem, which we prove in Section~\ref{section5}.
\begin{theom}Let $q\in \mathbb{N}$. The Hilbert series of $\mathbf{U}(\lie n^+)^{\otimes q}$ is given by
$$H(\mathbf{U}(\lie n^+)^{\otimes q}) = \sum\limits_{\bold k\in \mathbb{Z}_+^{|I|}}(-1)^{\sum\limits_{i\in I}k_i} \pi_{\bold k}^G(-q)\ e^{\eta(\bold k)}.$$
\end{theom}

The paper is organized as follows. In Section~\ref{section2} we define Borcherds algebras and state the denominator identity. 
In Section~\ref{section3} we define the generalized chromatic polynomial and give a combinatorial formula for the root multiplicities. In Section~\ref{section4} we construct a family of bases for the root spaces using the results of Section~\ref{section3}. In the last section we use the denominator identity to compute the Hilbert series of free partially commutative Lie algebras.

{\em Acknowledgements. This work was begun when the third author was visiting the University of Cologne, Germany and 
part of this work was done when the third author was at The Weizmann Institute of Science as a postdoctoral fellow.
The author is grateful to the University of Cologne and The Weizmann Institute of Science for the superb working environment.
The first and third authors thank Sankaran Viswanath for many helpful discussions.}

%%%%%%%%%%%%%%%%%%%%%%%%%%%%%%%%%%%%%%%%%%%%%%%%%
%%%%%%%%%%%%%%%%%%%%%%%%%%%%%%%%%%%%%%%%%%%%%%%%%
%%%%%%%%%%%%%%%%%%%%%%%%%%%%%%%%%%%%%%%%%%%%%%%%%
%%%%%%%%%%%%%%%%%%%%%%%%%%%%%%%%%%%%%%%%%%%%%%%%%
%%%%%%%%%%%%%%%%%%%%%%%%%%%%%%%%%%%%%%%%%%%%%%%%%

\section{Borcherds algebras and the denominator identity}\label{section2}
\subsection{} We denote the set of complex numbers by $\bc$ and, respectively, the set of integers, non--negative integers, and positive integers  by $\bz$, $\bz_+$, and $\bn$. Unless otherwise stated, all the vector spaces considered in this paper are $\bc$--vector spaces.
\subsection{}
We recall the definition of Borcherds algebras; also called generalized Kac--Moody algebras. For more details, we refer the reader to \cite{Bor88,Ej96,K90} and the references therein. A real matrix  
$A=(a_{ij})_{i,j\in I}$ indexed by a finite or countably infinite set, which we identify with $I=\{1, \dots, n\}$ or $\mathbb{Z}_+$, is said to be a \textit{Borcherds--Cartan\ matrix}
if the following conditions are satisfied for all $i,j\in I$:
\begin{enumerate}
\item $A$ is symmetrizable
\item $a_{ii}=2$ or $a_{ii}\leq 0$
\item $a_{ij}\leq 0$ if $i\neq j$ and $a_{ij}\in\mathbb{Z}$ if $a_{ii}=2$
\item $a_{ij}=0$ if and only if $a_{ji}=0$.
\end{enumerate} 
Recall that a matrix $A$ is called symmetrizable if there exists a diagonal matrix $D=\mathrm{diag}(\epsilon_i,i\in I)$ with positive entries such that $DA$ is symmetric. Set $I^\text{re}=\{i\in I: a_{ii}=2\}$ and $I^{\text{im}}=I\backslash I^\text{re}$. The Borcherds algebra $\lie g=\mathfrak{g}(A)$ associated to a Borcherds--Cartan matrix $A$ is the Lie algebra generated by $e_i, f_i, h_i$, $i\in I$ with the following defining relations:
\begin{enumerate}
 \item[(R1)] $[h_i, h_j]=0$ for $i,j\in I$
 \item[(R2)] $[h_i, e_k]=a_{i,k}e_i$,  $[h_i, f_k]=-a_{i,k}f_i$ for $i,k\in I$
 \item[(R3)] $[e_i, f_j]=\delta_{ij}h_i$ for $i, j\in I$
 \item[(R4)] $(\text{ad }e_i)^{1-a_{ij}}e_j=0$, $(\text{ad }f_i)^{1-a_{ij}}f_j=0$ if $i\in I^\text{re}$ and $i\neq j$
 \item[(R5)] $[e_i, e_j]=0$ and $[f_i, f_j]=0$ if $i, j \in I^\text{im}$ and $a_{ij}=0$.
\end{enumerate}
\begin{rem} 
If $i\in I$ is such that $a_{i,i}=0$, the subalgebra spanned by the elements $h_i,e_i,f_i$ is a Heisenberg algebra and otherwise this subalgebra is isomorphic to $\mathfrak{sl}_2$ (possibly after rescaling $e_i$ and $h_i$).
\end{rem}
\subsection{}We collect some elementary properties of Borcherds algebras; see \cite[Proposition 1.5]{Ej96} for more details and proofs. We have that $\lie g$ is $\mathbb{Z}^{I}$--graded by giving $h_i$ degree $(0,0,\dots)$, $e_i$ degree $(0,\dots,0,1,0,\dots)$ and $f_i$ degree $(0,\dots,0,-1,0,\dots)$ where $\pm 1$ appears at the $i$--th position. For a sequence $(n_1,n_2,\dots)$, we denote by $\lie g(n_1,n_2,\dots)$ the corresponding graded piece; note that $\lie g(n_1,n_2,\dots)=0$ unless finitely many of the $n_i$ are non--zero. Let $\mathfrak{h}$ be the abelian subalgebra spanned by the $h_i$, $i\in I$ and let $\mathfrak{E}$ the space of commuting derivations of $\lie g$ spanned by the $D_i$, $i\in I$, where $D_i$ denotes the derivation that acts on $\lie g(n_1,n_2,\dots)$ as multiplication by the scalar $n_i$. Note that the abelian subalgebra $\mathfrak{E}\ltimes \mathfrak{h}$ of $\mathfrak{E}\ltimes \mathfrak{g}$ acts by scalars on $\lie g(n_1,n_2,\dots)$ and we have a root space decomposition:
\begin{equation}\label{rootdec}\mathfrak{g}=\bigoplus _{\alpha \in (\mathfrak{E}\ltimes \mathfrak{h})^*}
\mathfrak{g}_{\alpha }, \ \mathrm{where} \ \mathfrak{g}_{\alpha }
:=\{ x\in \mathfrak{g}\ |\ [h, x]=\alpha(h) x \ 
\mathrm{for\ all}\ h\in \mathfrak{E}\ltimes \mathfrak{h} \}.\end{equation}
Define $\Pi=\{\alpha_i\}_{i\in I}\subset (\mathfrak{E}\ltimes\lie h)^{*}$ by $\alpha_j((D_k,h_i))=\delta_{k,j}+a_{i,j}$ and set
$$Q:=\bigoplus _{i\in I}\mathbb{Z}\alpha _i,\ \ Q_+ :=\sum _{i\in I}\mathbb{Z}_{+}\alpha _i.$$
Denote by $\Delta :=\{ \alpha \in (\mathfrak{E}\ltimes\lie h)^*\backslash \{0\} \mid \mathfrak{g}_{\alpha }\neq 0\}$ the set of roots, and by $\Delta_+$ the set of roots which are non--negative integral linear combinations of the $\alpha_i^{'}s$, called the positive roots. The elements in $\Pi$ are called the simple roots; we call $\Pi^\mathrm{re}:=\{\alpha_i: i\in I^\mathrm{re}\}$ the set of real simple roots and $\Pi^\mathrm{im}=\Pi\backslash \Pi^\mathrm{re}$ the set of imaginary simple roots. One of the important properties of Borcherds algebras is that $\Delta =\Delta_+ \sqcup - \Delta_+$ and 
$$\mathfrak{g}_0=\mathfrak{h},\ \ \lie g_\alpha=\lie g(n_1,n_2,\dots),\ \text{ if }\ \alpha=\sum_{i\in I} n_i\alpha_i\in \Delta.$$

Moreover, we have a triangular decomposition
$$\lie g\cong \lie n^{-}\oplus \lie h \oplus \lie n^+,$$
where $\lie n^{+}$ (resp. $\lie n^{-}$) is the free Lie algebra generated by $e_i,\ i\in I$ (resp. $f_i,\ i\in I$) with defining relations 
$$(\text{ad }e_i)^{1-a_{ij}}e_j=0\ (\text{resp. } (\text{ad }f_i)^{1-a_{ij}}f_j=0) \text{ for } i\in I^\text{re} \text{ and } i\neq j$$
and
$$[e_i, e_j]=0\ (\text{resp. } [f_i, f_j]=0) \text{ for } i, j \in I^\text{im} \text{ and } a_{ij}=0.$$

In view of \eqref{rootdec} we have
$$\lie n^{\pm}=\bigoplus_{\alpha \in \pm\Delta_{+}}
\mathfrak{g}_{\alpha}.$$
Finally, given $\gamma=\sum_{i\in I}n_i\alpha_i\in Q_+$ (only finitely many $n_i$ are non--zero), we set $\text{ht}(\gamma)=\sum\limits_{i\in I}n_i.$ 
The following lemma has been proved in \cite[Corollary 11.13.1]{K90} for a finite index set $I$ and remains valid without any modification for countable $I$. We will need this result in Section~\ref{section3}.
\begin{lem}\label{rootspro}
Let $i\in I^{\text{im}}$ and $\alpha\in \Delta_+\backslash\{\alpha_i\}$ such that $\alpha(h_i)<0$. Then $\alpha+j\alpha_i\in \Delta_+$ for all $j\in \bz_+$.
\qed\hfill
\end{lem}
 
\begin{rem} 
Although Kac--Moody algebras are constructed similarly as Borcherds algebras using generalized Cartan matrices (see \cite{K90} for details), the theory of Borcherds algebras includes examples which behave in a very different way from Kac--Moody algebras. The main point of difference is that one is allowed to have imaginary simple roots. 
\end{rem}

\subsection{}
We denote by $R=Q\otimes_{\bz} \mathbb{R}$ the real vector space spanned by $\Delta$. There exists a symmetric bilinear form on $R$ given by $(\alpha_i, \alpha_j)=\epsilon_i a_{ij}$ for $i, j\in I.$ For $i\in I^{\text{re}}$, define the linear isomorphism $\bold{s}_i$ of $R$ by 
$$\bold{s}_i(\lambda)=\lambda-\lambda(h_i)\alpha_i=\lambda-2\frac{(\lambda,\alpha_i)}{(\alpha_i,\alpha_i)}\alpha_i,\ \ \lambda\in R.$$ The Weyl group $W$ of $\mathfrak{g}$ is the subgroup of $\mathrm{GL}(R)$ generated by the simple reflections $\bold{s}_i$, $i\in I^\mathrm{re}$. Note that $W$ is a Coxeter group with canonical generators $\bold{s}_i, i\in I^\mathrm{re}$ and the above bilinear form is $W$--invariant. We denote by $\ell(w)=\mathrm{min}\{k\in \mathbb{N}: w=\bold{s}_{i_1}\cdots \bold{s}_{i_k}\}$ the length of $w\in W$ and $w=\bold{s}_{i_1}\cdots \bold{s}_{i_k}$ with $k=\ell(w)$ is called a reduced expression. We denote by $\Delta^\mathrm{re}=W(\Pi^\mathrm{re})$ the set of real roots and $\Delta^\mathrm{im}=\Delta\backslash \Delta^\mathrm{re}$ the set of imaginary roots. Equivalently, a root $\alpha$ is imaginary if and only if $(\alpha, \alpha)\le 0$ and else real. We can extend $(.,.)$ to a symmetric form on $(\mathfrak{E}\ltimes \lie h)^*$ satisfying $(\lambda,\alpha_i)=\lambda(\epsilon_ih_i)$ and also $\bold{s}_i$ to a linear isomorphism of $(\mathfrak{E}\ltimes \lie h)^*$ by 
$$\bold{s}_i(\lambda)=\lambda-\lambda(h_i)\alpha_i,\ \ \lambda\in (\mathfrak{E}\ltimes \lie h)^*.$$
Let $\rho$ be any element of $(\mathfrak{E}\ltimes \lie h)^*$ satisfying $2(\rho,\alpha_i)=(\alpha_i,\alpha_i)$ for all $i\in I$.
 The following denominator identity has been proved in \cite{Bor88}, see also \cite[Theorem 3.16]{Ej96}.
\subsection*{The denominator identity} 
\begin{eqnarray}\label{denominator}
U:=  \sum_{w \in W } (-1)^{\ell(w)} \sum_{ \gamma \in
\Omega} (-1)^{\mathrm{ht}(\gamma)} e^{w(\rho -\gamma)-\rho } & =& \prod_{\alpha \in \Delta_+} (1 - e^{-\alpha})^{\dim \mathfrak{g}_{\alpha}} 
\label{eq:denom} \end{eqnarray}  
where $\Omega$ is the set of all $\gamma\in Q_+$ such that $\gamma$ is a finite sum of mutually orthogonal distinct imaginary simple roots. Note that $0\in \Omega$ and $\alpha\in \Omega,$ if $\alpha$ is imaginary. 
%%%%%%%%%%%%%%%%%%%%%%%%%%%%%%%%%%%%%%%%%%%%%%%%%
%%%%%%%%%%%%%%%%%%%%%%%%%%%%%%%%%%%%%%%%%%%%%%%%%
%%%%%%%%%%%%%%%%%%%%%%%%%%%%%%%%%%%%%%%%%%%%%%%%%
%%%%%%%%%%%%%%%%%%%%%%%%%%%%%%%%%%%%%%%%%%%%%%%%%
\section{Proper multicoloring, root multiplicities and chromatic polynomials}\label{section3}
In this section we associate a graph $G$ to a given Borcherds algebra $\lie g$ and establish a connection between the root multiplicities and the generalized chromatic polynomial of $G$. The main results of this section are Theorem~\ref{mainthmch} and the closed formula stated in Corollary~\ref{recursionmult}.
\subsection{}For a given Borcherds algebra $\lie g$ we associate a graph $G$ as follows: $G$ has vertex set $I$, with an edge between two vertices $i$ and $j$ iff 
$a_{ij}\neq 0$ for $i,j\in I$, $i\neq j$.
Note that $G$ is a simple (finite or infinite) graph which we call the $\mathrm{graph \ of} \ \lie g$.
The edge set of $G$ is denoted by $E(G)$ and $e(i, j)$ denotes the edge between the nodes $i$ and $j$. 
A finite subset $S\subset I$ is said to be connected if the corresponding subgraph generated by $S$ is connected. 
\subsection{}\label{multicoloring}In this subsection we discuss the notion of a proper vertex multicoloring of a given graph, 
which is a generalization of the well--known graph coloring. For more details about multicoloring of a graph we refer to \cite{HMK04}.
For any finite set $S$, let $\mathcal{P}(S)$ be the power set of $S$.

\begin{defn}Let $\mathbf{k}=(k_i: i\in I)$ a tuple of non--negative integers such that $|\{i\in I: k_i\neq 0\}|<\infty$. 
We call a map $\tau: I\rightarrow \mathcal{P}\big(\{1,\dots,q\}\big)$ a proper vertex multicoloring of $G$ associated to $\bold k$ if the following conditions are satisfied:
\begin{enumerate}
\item[(i)]  For all $i\in I$ we have $|\tau(i)|=k_i$ 
\item[(ii)] For all $i,j\in I$ such that $(i,j)\in E(G)$ we have $\tau(i)\cap \tau(j)=\emptyset$,
\end{enumerate}
where we understand $|\emptyset|=0$.
\end{defn}
It means that no two adjacent vertices share the same color. 
Note that if $S\subset I$, then the choice $k_i=1$ for $i\in S$ and $k_i=0$ otherwise corresponds to an ordinary graph coloring of the subgraph spanned by $S$. 
\begin{example}
We consider the graph 

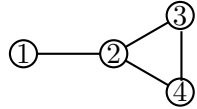
\hspace{3,4cm}
\begin{tikzpicture}[scale=.6]
    \draw (-1,0) node[anchor=east]{};
    \draw[xshift=3 cm,thick] (3 cm,0) circle(.3cm) node[]{1};
    \draw[xshift=4 cm,thick] (4 cm,0) circle(.3cm)node[]{2};
    \draw[xshift=8 cm,thick] (30: 17 mm) circle (.3cm)node[]{3};
    \draw[xshift=8 cm,thick] (-30: 17 mm) circle (.3cm)node[]{4};
  %  \draw[dotted,thick] (0.3 cm,0) -- +(1.4 cm,0);
    \foreach \y in {3.15,...,3.15}
    \draw[xshift=\y cm,thick] (\y cm,0) -- +(1.4 cm,0);
    \draw[xshift=8 cm,thick] (30: 3 mm) -- (30: 14 mm);
    \draw[xshift=8 cm,thick] (-30: 3 mm) -- (-30: 14 mm);
      \draw[thick] (95mm, -6mm) -- +(0, 11mm);
\end{tikzpicture}

We allow 3 colors, say blue, red and green and fix $\bold k=(2,1,1,1)$. Below we have listed all proper vertex multicolorings

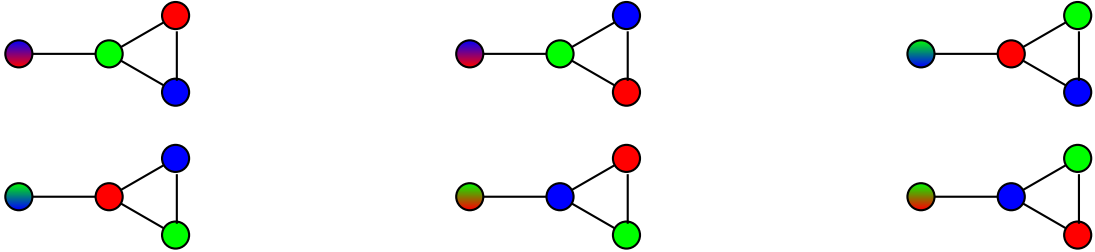
\hspace{-4cm}
\begin{tikzpicture}[scale=.6]
    \draw (-1,0) node[anchor=east]{};
    \draw[top color=blue,bottom color=red, xshift=3 cm,thick] (3 cm,0) circle(.3cm);
    \draw[fill=green,xshift=4 cm,thick] (4 cm,0) circle(.3cm);
    \draw[fill=red, xshift=8 cm,thick] (30: 17 mm) circle (.3cm);
    \draw[fill=blue, xshift=8 cm,thick] (-30: 17 mm) circle (.3cm);
  %  \draw[dotted,thick] (0.3 cm,0) -- +(1.4 cm,0);
    \foreach \y in {3.15,...,3.15}
    \draw[xshift=\y cm,thick] (\y cm,0) -- +(1.4 cm,0);
    \draw[xshift=8 cm,thick] (30: 3 mm) -- (30: 14 mm);
    \draw[xshift=8 cm,thick] (-30: 3 mm) -- (-30: 14 mm);
      \draw[thick] (95mm, -6mm) -- +(0, 11mm);
\end{tikzpicture}
\hspace{-1cm}
\begin{tikzpicture}[scale=.6]
    \draw (-1,0) node[anchor=east]{};
    \draw[top color=blue,bottom color=red, xshift=3 cm,thick] (3 cm,0) circle(.3cm);
    \draw[fill=green, xshift=4 cm,thick] (4 cm,0) circle(.3cm);
    \draw[fill=blue, xshift=8 cm,thick] (30: 17 mm) circle (.3cm);
    \draw[fill=red, xshift=8 cm,thick] (-30: 17 mm) circle (.3cm);
  %  \draw[dotted,thick] (0.3 cm,0) -- +(1.4 cm,0);
    \foreach \y in {3.15,...,3.15}
    \draw[xshift=\y cm,thick] (\y cm,0) -- +(1.4 cm,0);
    \draw[xshift=8 cm,thick] (30: 3 mm) -- (30: 14 mm);
    \draw[xshift=8 cm,thick] (-30: 3 mm) -- (-30: 14 mm);
      \draw[thick] (95mm, -6mm) -- +(0, 11mm);
\end{tikzpicture}
\hspace{-1cm}
\begin{tikzpicture}[scale=.6]
    \draw (-1,0) node[anchor=east]{};
    \draw[top color=green,bottom color=blue,xshift=3 cm,thick] (3 cm,0) circle(.3cm);
    \draw[fill=red, xshift=4 cm,thick] (4 cm,0) circle(.3cm);
    \draw[fill=green, xshift=8 cm,thick] (30: 17 mm) circle (.3cm);
    \draw[fill=blue, xshift=8 cm,thick] (-30: 17 mm) circle (.3cm);
  %  \draw[dotted,thick] (0.3 cm,0) -- +(1.4 cm,0);
    \foreach \y in {3.15,...,3.15}
    \draw[xshift=\y cm,thick] (\y cm,0) -- +(1.4 cm,0);
    \draw[xshift=8 cm,thick] (30: 3 mm) -- (30: 14 mm);
    \draw[xshift=8 cm,thick] (-30: 3 mm) -- (-30: 14 mm);
      \draw[thick] (95mm, -6mm) -- +(0, 11mm);
\end{tikzpicture}
\\

\hspace{-4cm}
\begin{tikzpicture}[scale=.6]
    \draw (-1,0) node[anchor=east]{};
    \draw[top color=green,bottom color=blue, xshift=3 cm,thick] (3 cm,0) circle(.3cm);
    \draw[fill=red, xshift=4 cm,thick] (4 cm,0) circle(.3cm);
    \draw[fill=blue, xshift=8 cm,thick] (30: 17 mm) circle (.3cm);
    \draw[fill=green, xshift=8 cm,thick] (-30: 17 mm) circle (.3cm);
  %  \draw[dotted,thick] (0.3 cm,0) -- +(1.4 cm,0);
    \foreach \y in {3.15,...,3.15}
    \draw[xshift=\y cm,thick] (\y cm,0) -- +(1.4 cm,0);
    \draw[xshift=8 cm,thick] (30: 3 mm) -- (30: 14 mm);
    \draw[xshift=8 cm,thick] (-30: 3 mm) -- (-30: 14 mm);
      \draw[thick] (95mm, -6mm) -- +(0, 11mm);
\end{tikzpicture}
\hspace{-1cm}
\begin{tikzpicture}[scale=.6]
    \draw (-1,0) node[anchor=east]{};
    \draw[top color=green,bottom color=red, xshift=3 cm,thick] (3 cm,0) circle(.3cm);
    \draw[fill=blue, xshift=4 cm,thick] (4 cm,0) circle(.3cm);
    \draw[fill=red, xshift=8 cm,thick] (30: 17 mm) circle (.3cm);
    \draw[fill=green, xshift=8 cm,thick] (-30: 17 mm) circle (.3cm);
  %  \draw[dotted,thick] (0.3 cm,0) -- +(1.4 cm,0);
    \foreach \y in {3.15,...,3.15}
    \draw[xshift=\y cm,thick] (\y cm,0) -- +(1.4 cm,0);
    \draw[xshift=8 cm,thick] (30: 3 mm) -- (30: 14 mm);
    \draw[xshift=8 cm,thick] (-30: 3 mm) -- (-30: 14 mm);
      \draw[thick] (95mm, -6mm) -- +(0, 11mm);
\end{tikzpicture}
\hspace{-1cm}
\begin{tikzpicture}[scale=.6]
    \draw (-1,0) node[anchor=east]{};
    \draw[top color=green,bottom color=red, xshift=3 cm,thick] (3 cm,0) circle(.3cm);
    \draw[fill=blue, xshift=4 cm,thick] (4 cm,0) circle(.3cm);
    \draw[fill=green,xshift=8 cm,thick] (30: 17 mm) circle (.3cm);
    \draw[fill=red, xshift=8 cm,thick] (-30: 17 mm) circle (.3cm);
  %  \draw[dotted,thick] (0.3 cm,0) -- +(1.4 cm,0);
    \foreach \y in {3.15,...,3.15}
    \draw[xshift=\y cm,thick] (\y cm,0) -- +(1.4 cm,0);
    \draw[xshift=8 cm,thick] (30: 3 mm) -- (30: 14 mm);
    \draw[xshift=8 cm,thick] (-30: 3 mm) -- (-30: 14 mm);
      \draw[thick] (95mm, -6mm) -- +(0, 11mm);
\end{tikzpicture}
\end{example}
Multicoloring play an important role in algebraic graph theory. One of the favorite examples 
is scheduling dependent jobs on multiple machines. When all jobs have the same execution times, this is modeled by a graph coloring problem and as a graph multicoloring problem for arbitrary execution times. The vertices in the graph represent the jobs and an edge in the graph between two vertices forbids scheduling these jobs simultaneously. For more details and examples we refer to \cite{HMK04}. 

\subsection{}\label{chromatic}
The number of ways a graph $G$ can be $\bold k$--multicolored using $q$ colors is a polynomial in $q$, called the generalized chromatic polynomial $\pi_\mathbf{k}^G(q)$. 
The generalized chromatic polynomial has the following well--known description. 
We denote by $P_k(\bold k,G)$ the set of all ordered $k$--tuples $(P_1,\dots,P_k)$ such that:
\begin{enumerate}
 \item[(i)] each $P_i$ is a non--empty independent subset of $I$, i.e. no two vertices have an edge between them; and

 \item[(ii)] the disjoint union of  $P_1,\cdots, P_k$ is equal to the multiset $\{\underbrace{\alpha_i,\dots, \alpha_i}_{k_i\, \mathrm{times}}: i\in I \}$.
\end{enumerate} 
Then we have
\begin{equation}\label{defgenchr}\pi^G_\mathbf{k}(q)= \sum\limits_{k\ge0}|P_k(\bold k, G)| \, {q \choose k}.\end{equation}

The definition of the generalized chromatic polynomial depends on $q\in \mathbb N$; however it makes sense to consider the evaluation at any integer. 
A famous result of Stanley \cite{S73} says that the evaluation at $q=-1$ of the chromatic polynomial counts (up to a sign) the number of acyclic orientations. 
In Section~\ref{section5} we give a Lie theoretic proof of Stanley's result and more generally we prove his reciprocity theorem of chromatic polynomials.
\subsection{}
\textit{Fix for the rest of this paper a tuple of non--negative integers $\mathbf{k}=(k_i: i\in I)$ such that $k_i\leq 1$ for $i\in I^{\text{re}}$ 
and $1<|\{i\in I: k_i\neq 0\}|<\infty$.} We set $\eta(\mathbf{k})=\sum_{i\in I}k_i\alpha_i\in Q_+$.
\begin{defn}\label{ing}Let $L_{G}(\bold k)$ be the weighted bond lattice of $G$, which is the set of 
$\mathbf{J}=\{J_1,\dots,J_k\}$  satisfying the following properties:
\begin{enumerate}
\item[(i)] $\bold J$ is a multiset, i.e. we allow $J_i=J_j$ for $i\neq j$

\item[(ii)] each $J_i$ is a multiset and the subgraph spanned by the underlying set of $J_i$ is connected subgraph of $G$ for each $1\le i\le k$ and

 \item[(iii)] the disjoint union of $J_1\dot{\cup} \cdots \dot{\cup} J_k=\{\underbrace{\alpha_i,\dots, \alpha_i}_{k_i\, \mathrm{times}}: i\in I\}$.
 
\end{enumerate}

 For $\mathbf{J}\in L_{G}(\bold k)$ we denote by $D(J_i,\mathbf{J})$ the multiplicity of $J_i$ in $\mathbf{J}$ and set $\text{mult}(\beta(J_i))=\text{dim }\mathfrak{g}_{\beta(J_i)}$, where $\beta(J_i)=\sum_{\alpha\in J_i} \alpha$.
\end{defn}
 We record the following lemma which will be needed later.
\begin{lem}\label{bij}Let $\mathcal{P}$ be the collection of sets $\gamma=\{\beta_1,\dots,\beta_r\}$ (we allow $\beta_i=\beta_j$ for $i\neq j$) such that each $\beta_i\in\Delta_+$ and $\beta_1+\dots+\beta_r=\eta(\bold k)$.
The map $\Psi: L_{G}(\bold k)\rightarrow \mathcal{P}$ defined by $\{J_1,\dots,J_k\}\mapsto \{\beta(J_1),\dots,\beta(J_k)\}$ is a bijection.
\begin{pf}
If $\alpha\in \big(Q_+\cap\sum_{j\in \Pi^{\mathrm{re}}} \mathbb{Z}_{\leq 1}\alpha_j\big)$ is non--zero and the support of $\alpha$ is connected, then $\alpha\in\Delta^{\mathrm{re}}_+$. Moreover, if $\alpha\in \Delta_+$ and $\alpha_i\in \Pi^{\mathrm{im}}$ is such that the support of $\alpha+\alpha_i$ is connected, then by Lemma~\ref{rootspro} we have that $\alpha+\alpha_i\in \Delta_+$. This shows that each $\beta(J_r)$ is a positive root and hence the map is well--defined. The map is obviously injective and since $\alpha\in\Delta_+$ implies that $\alpha$ has connected support, we also obtain that $\Psi$ is surjective.
\end{pf}
\end{lem}
\subsection{}
The rest of this section is dedicated to the proof of the following theorem, which relates the generalized chromatic polynomial with root multiplicities of Borcherds algebras. Hence we obtain a Lie theoretic interpretation of generalized chromatic polynomials.
\begin{thm}\label{mainthmch}
Let $G$ be the graph of a Borcherds algebra. Then
$$
\pi^G_{\bold k}(q)=\sum_{\mathbf{J}\in L_{G}(\bold k)} (-1)^{\mathrm{ht}(\eta(\bold k))+|\mathbf{J}|}\prod_{J\in\bold J}\binom{q\text{ mult}(\beta(J))}{D(J,\mathbf{J})}.
$$
\end{thm}
\begin{rem}
The above theorem is a generalization of \cite[Theorem 1.1]{VV15} where the authors considered the special case when $\lie g$ is a Kac--Moody algebra and $k_i=1$ for all $i\in I$. 
\end{rem}
\subsection{}
For a Weyl group element $w\in W$,
we fix a reduced word $w=\bold {s}_{i_1}\cdots \bold{s}_{i_k}$ and let $I(w)=\{\alpha_{i_1},\dots,\alpha_{i_k}\}$. 
Note that $I(w)$ is independent of the choice of the reduced expression of $w$. For $\gamma\in \Omega$ we set $I(\gamma)=\{\alpha\in \Pi^{\mathrm im} : \mbox{ $\alpha$ is a summand of $\gamma$}\}$ and 
$$\mathcal{J}(\gamma)=\{ w\in W\backslash \{e\}: I(w)\cup I(\gamma) \mbox{ is an independent set}\}.$$
 Note that $\mathcal{J}(0)$ gives the set of independent subsets of $\Pi^\mathrm{re}$. The following lemma is the generalization of \cite[Lemma 2.3]{VV15} in the setting of Borcherds algebras.
\begin{lem}\label{helplem}
Let $w\in W$ and $\gamma\in \Omega$. We write $\rho-w(\rho)+w(\gamma)=\sum_{\alpha\in\Pi} b_{\alpha}(w,\gamma)\alpha$. Then we have

\begin{enumerate}
\item[(i)] $b_{\alpha}(w,\gamma)\in \mathbb{Z}_{+}$ for all $\alpha\in \Pi$ and $b_{\alpha}(w,\gamma)= 0$ if $\alpha\notin I(w)\cup I(\gamma)$,
\item[(ii)] $I(w)=\{\alpha\in \Pi^{\mathrm re} : b_{\alpha}(w,\gamma)\geq 1\}$ and  $b_{\alpha}(w,\gamma)=1$ if $\alpha\in I(\gamma)$,
\item[(iii)] If $w\in \mathcal{J}(\gamma)$, then $b_{\alpha}(w,\gamma)=1$ for all $\alpha\in I(w)\cup I(\gamma)$ and $b_{\alpha}(w,\gamma)=0$ else,
\item[(iv)] If $w\notin \mathcal{J}(\gamma)\cup \{e\}$, then there exists $\alpha\in \Pi^{\mathrm re}$ such that $b_{\alpha}(w,\gamma)>1$.

\end{enumerate}
\begin{pf}
We start proving (i) and (ii) by induction on $\ell(w)$. If $\ell(w)=0$, the statement is obvious. So let $\alpha\in \Pi^\mathrm{re}$ such that $w=s_{\alpha}u$ and $\ell(w)=\ell(u)+1$. Then
\begin{align}\notag \rho-w(\rho)+w(\gamma)&=\rho-s_{\alpha}u(\rho)+s_{\alpha}u(\gamma)&\\&\label{1}=\rho-u(\rho)+u(\gamma)+2\frac{(\rho,u^{-1}\alpha)}{(\alpha,\alpha)}\alpha-2\frac{(\gamma,u^{-1}\alpha)}{(\alpha,\alpha)}\alpha.
\end{align}
So by our induction hypothesis we know $\rho-u(\rho)+u(\gamma)$ has the required property and since $\ell(w)=\ell(u)+1$, we also know  $u^{-1}\alpha\in \Delta^{\mathrm re} \cap \Delta_+$. 
Note that $(\rho,\alpha_i)=\frac{1}{2}(\alpha_i,\alpha_i)$ for all $\alpha_i\in \Pi^{\text{re}}$ implies  $2\frac{(\rho,u^{-1}\alpha)}{(\alpha,\alpha)}\in \mathbb{N}$. 
Furthermore, $\gamma$ is a sum of imaginary simple roots and $a_{ij}\leq 0$ whenever $i\neq j$. 
Hence $-2\frac{(\gamma,u^{-1}\alpha)}{(\alpha,\alpha)}\in \mathbb{Z}_{+}$ and the proof of (i) and (ii) is done, since $I(w)=I(u)\cup \{\alpha\}$.
If $w\in \mathcal{J}(\gamma)$ and $\alpha\in I(w)\cup I(\gamma)$, we have $(\rho,u^{-1}\alpha)=(\rho,\alpha)=\frac{1}{2}(\alpha,\alpha)$ and $(\gamma,u^{-1}\alpha)=(\gamma,\alpha)=0$. 
So part (iii) follows from \eqref{1} and an induction argument on $\ell(w)$ since $I(w)=I(u)\cup\{\alpha\}$ and $u\in  \mathcal{J}(\gamma)$.
It remains to prove part (iv), which will be proved again by induction. If $w=s_{\alpha}$ we have
$$\rho-w(\rho)+w(\gamma)=\alpha+\gamma-2\frac{(\gamma,\alpha)}{(\alpha,\alpha)}\alpha.$$
Since $w\notin \mathcal{J}(\gamma)\cup \{e\}$ we get $-2\frac{(\gamma,\alpha)}{(\alpha,\alpha)}\in\mathbb{N}$. 
For the induction step we write $w=s_{\alpha}u$. We have $s_\alpha\notin \mathcal{J}(\gamma)\cup \{e\}$ or $u\notin\mathcal{J}(\gamma)\cup \{e\}$. 
In the latter case we are done by using the induction hypothesis and \eqref{1}. 
Otherwise we can assume that $u\in \mathcal{J}(\gamma)\cup \{e\}$ and hence $(\gamma,u^{-1}\alpha)=(\gamma,\alpha)<0$ and $2\frac{(\rho,u^{-1}\alpha)}{(\alpha,\alpha)}\in \mathbb{N}$. Now interpreting this in \eqref{1} gives the result.
\end{pf}
\end{lem}
\subsection{}The following proposition is an easy consequence of Lemma~\ref{helplem} and will be needed in the proof of Theorem~\ref{mainthmch}.
Recall that $U$ is the sumside of the denominator identity \ref{denominator}.
\begin{prop}\label{helprop}
Let $q\in \mathbb{Z}$. We have
 
$$U^q[e^{-\eta(\bold k)}]= (-1)^{\mathrm{ht}(\eta(\bold k))}\  \pi^G_\mathbf{k}(q),$$
where $U^q[e^{-\eta(\bold k)}]$ denotes the coefficient of $e^{-\eta(\bold k)}$ in $U^q$.
\begin{pf} If $q=0$, then there is nothing to prove. So assume that $0\neq q\in \mathbb{Z}.$
 We have $$U^q=\sum\limits_{k\ge 0}{q \choose k}\, (U-1)^k, \ \ \text{where ${q\choose k} = \frac{q(q-1)\cdots (q-(k-1)))}{k!}$}.$$ 
 Note that $${-q\choose k} = 
 (-1)^k{q+k-1\choose k},\ \ \text{for $q\in \mathbb{N}$}.$$
From Lemma~\ref{helplem} we get 
$$w(\rho)-\rho-w(\gamma)=-\gamma-\sum_{\alpha\in I(w)}\alpha,\ \ \text{for $w\in \mathcal{J}(\gamma)\cup\{e\}$}$$ and thus $(U-1)^k[e^{-\eta(\bold k)}]$ is equal to 
$$\Bigg(\sum_{w \in \mathcal{J}(0)} (-1)^{\ell(w)}e^{-\sum_{\alpha\in I(w)}\alpha} + \sum_{ \gamma \in
\Omega \backslash \{0\} }(-1)^{\mathrm{ht}(\gamma)}\sum_{w \in \mathcal{J}(\gamma)\cup\{e\} } (-1)^{\ell(w)}  e^{-\gamma-\sum_{\alpha\in I(w)}\alpha}\Bigg)^k[e^{-\eta(\bold k)}].$$
Hence the coefficient is given by
$$\sum_{\substack{(\gamma_1,\dots,\gamma_k)\\(w_1,\dots,w_k)}}(-1)^{\sum_{i=1}^k\mathrm{ht}(\gamma_i)}(-1)^{\ell(w_1\cdots w_k)},$$
where the sum ranges over all $k$--tuples $(\gamma_1,\dots,\gamma_k)\in \Omega^k$ (repetition is allowed) and $(w_1,\dots,w_k)$ such that
\begin{align*}
&\bullet \ w_i\in \mathcal{J}(\gamma_i)\cup \{e\},\text{  $1\le i\le k$},&\\&
\bullet \ I(w_1)\ \dot{\cup} \cdots \dot{\cup}\ I(w_k)=\{\alpha_i: i\in I^{\mathrm{re}}, k_i=1\},&\\&
\bullet \ I(w_i)\cup I(\gamma_i)\neq \emptyset \ \text{ for each $1\le i\le k$},&\\&
\bullet \ \gamma_1+\cdots+\gamma_k=\sum_{i\in I^{\mathrm{im}}}k_i\alpha_i. 
\end{align*}

It follows that $\big(I(w_1)\cup I(\gamma_1),\dots,I(w_k)\cup I(\gamma_k)\big)\in P_k\big(\bold k,G\big)$ and each element is obtained in this way. So the sum ranges over all elements in $P_k(\bold k,G).$
Since $w_1\cdots w_k$ is a subword of a Coxeter element we get
$$(-1)^{\ell(w_1\cdots w_k)}=(-1)^{|\{i\in I^{\mathrm{re}} : k_i=1\}|},$$
and hence $(U-1)^k[e^{-\eta(\bold k)}]$ is equal to $(-1)^{\mathrm{ht}(\eta(\bold k))} |P_k(\bold k, G)|$ which finishes the proof.
\end{pf}
\end{prop}

%%%%%%%%%%%%%%%%%%%%%%%%%%%%%%%%%%%%%%%%%%%%%%%%%%%%%%%%%%%%%%%%%%%%%%%%%%%%%%%%%%%%%%%%%%%%%%%%%%%%%%%%%%%%%%%%%%%%%%%%%%%%%%%%%%%%%%%%%%%%%%%%%%%%%%%%%%%%%%%%%%%%%%%%
\subsection{}
Now we are able to prove Theorem~\ref{mainthmch} by using the denominator identity \eqref{eq:denom}. Proposition~\ref{helprop} and \eqref{eq:denom} together imply that the generalized chromatic polynomial $\pi^G_\mathbf{k}(q)$ is given by the coefficient of $e^{-\eta(\bold k)}$ in 
\begin{equation}\label{expa}(-1)^{\mathrm{ht}(\eta(\bold k))}\prod_{\alpha\in \Delta_+}(1-e^{-\alpha})^{q\text{ dim }\mathfrak{g}_{\alpha}}.\end{equation}
Expanding \eqref{expa} and using Lemma~\ref{bij} finishes the proof.

%%%%%%%%%%%%%%%%%%%%%%%%%%%%%%%%%%%%%%%%%%%%%%%%%%%%%%%%%%%%%%%%%%%%%%%%%%%%%%%%%%%%%%%%%%%%%%%%%%%%%%%%%%%%%%%%%%%%%%%%%%%%%%%%%%%%%%%%%%%%%%%%%%%%%%%%%%%%%%%%%%%%%%%%

\subsection{}
In this subsection we prove a corollary which gives a combinatorial formula for certain root multiplicities. We consider the algebra of formal power series $\mathcal{A}:=\mathbb{C}[[X_i : i\in I]]$. For a formal power series $\zeta\in\mathcal{A}$ with constant term 1, its logarithm $\text{log}(\zeta)=-\sum_{k\geq 1}\frac{(1-\zeta)^k}{k}$ is well--defined.

\begin{cor}\label{recursionmult}
We have
\begin{equation}\label{mult}
  \text{ mult } \eta(\bold k) = \sum\limits_{\ell | \bold k}\frac{\mu(\ell)}{\ell}\ |\pi^G_{\bold k/\ell}(q)[q]|,\end{equation}
where $|\pi^G_{\bold k}(q)[q]|$ denotes the absolute value of the coefficient of $q$ in $\pi^G_{\bold k}(q)$ and $\mu$ is the M\"{o}bius function.

\begin{pf}
We consider $U$ as an element of $\mathbb{C}[[e^{-\alpha_i} : i\in I]]$. 
From the proof of Proposition~\ref{helprop} we obtain that the coefficient of $e^{-\eta(\bold k)}$ in $-\text{log } U$ equals
$$(-1)^{\mathrm{ht}(\eta(\bold k))}\sum\limits_{k\ge 1}\frac{(-1)^k}{k}|P_k(\bold k,G)|$$ which by definition \eqref{defgenchr} is equal to $|\pi^G_{\bold k}(q)[q]|$.
Now applying $-\text{log }$ to the right hand side of the denominator identity \eqref{eq:denom} gives 
\begin{equation}\label{recursss}\sum_{\substack{\ell\in \mathbb{N}\\ \ell | \bold k}} \frac{1}{\ell}\text{ mult } \eta\left(\bold k/\ell\right)=
|\pi^G_{\bold k}(q)[q]|.\end{equation} The statement of the corollary is now an easy consequence of the M\"{o}bius inversion formula.
\end{pf}
\end{cor}

\begin{rem}
The generalized chromatic polynomial $\pi^G_{\bold k}(q)$ can be computed explicitly for many families of graphs and hence \eqref{mult} gives an effective method to compute the root multiplicities; the case of complete graphs and trees is treated in the end of this section. For chordal graphs we refer the reader to \cite{AG03}.
\end{rem} 

\subsection{}\label{examples}
We calculate the generalized chromatic polynomial for special families of graphs.
\begin{example}\mbox{}
\begin{enumerate}
\item Let $G=K_n$ be a complete graph with $n$ vertices and $\bold k=(k_1, \dots, k_n)$ be a tuple of positive integers. Vertex $1$ can receive any $k_1$ colors from the given $q$ colors. Vertex $2$ can not receive those $k_1$ colors that were assigned to vertex $1$ and that is the only restriction that we have. Hence vertex 2
  can receive any $k_2$ colors from the remaining $q-k_1$ colors. Similarly the vertex $3$ can receive any $k_3$ colors from the remaining $q-(k_1+k_2)$ colors. Continuing in this way, we get that the generalized chromatic polynomial is given by
 $$\pi^{K_n}_{\bold k}(q)= {q\choose k_1}{q-k_{1}\choose k_2}{q-(k_{1}+k_2)\choose k_3}\cdots {q-(k_1+\cdots+k_{n-1})\choose k_n}.$$
In particular $\pi^{K_n}_{\bold k}(q)[q]=\frac{(k_1+\cdots k_n-1)!}{k_1!\cdots k_n!}.$
Hence we recover Witt's formula proved in \cite{Wi37}:
$$\text{ mult } \eta(\bold k) = \frac{1}{\rm{ht }\bold k}\sum\limits_{\ell | \bold k}\mu(\ell)\frac{(\rm{ht }\bold k/\ell)!}{(\bold k/\ell)!}.$$
 
\item Let $G=T_n$ be a tree with $n$ vertices and $\bold k=(k_1,\dots,k_n)$ a tuple of positive integers. Assume that the vertex set $I=\{1,\dots,n\}$ of $G$ is ordered in such a way that the vertex $i$ is a leaf of the subgraph of $G$ spanned by the vertices $\{i, i+1,\dots,n\}$. Further, let $i'$ the unique vertex adjacent to $i$ for each $1\le i\le n-1.$  We denote by $G'$ the subgraph obtained from $G$ by deleting vertex $1$ and claim that each vertex multicoloring of $G'$ gives ${q-k_{1'}\choose k_1}$ distinct vertex multicolorings of $G$. Fix a multicoloring of $G'$, then vertex $1$ is colored with $k_{1}$ distinct multicolors. Now
 it is easy to see that vertex $1$ can not be colored by the colors which are used to color vertex $1'$ and this is the only restriction that we have. Hence we can choose any $k_1$ colors among the $q-k_{1'}$ remaining colors to color vertex $1$. This proves that
 $$\pi^{T_n}_{\bold k}(q)= {q-k_{1'}\choose k_1}\pi^{G'}_{\bold k'}(q),$$
 where $\bold k'=(k_2, \dots, k_n).$  Now we can repeat this procedure and obtain together with $\pi^{\{n\}}_{k_n}(q)={q\choose k_n}$ that
 $$\pi^{T_n}_{\bold k}(q)= {q-k_{1'}\choose k_1}{q-k_{2'}\choose k_2}\cdots {q-k_{(n-1)'}\choose k_{n-1}}{q\choose k_n}.$$
In particular we get a Witt type formula for trees:
 $$  \text{ mult } \eta(\bold k) = \sum\limits_{\ell | \bold k}\mu(\ell) {\frac{k_1+k_{1'}}{\ell}-1\choose \frac{k_1}{\ell}}{\frac{k_2+k_{2'}}{\ell}-1\choose \frac{k_2}{\ell}}\cdots
 {\frac{k_{n-1}+k_{(n-1)'}}{\ell}-1\choose \frac{k_{n-1}}{\ell}}(1/k_n).$$

\end{enumerate}

\end{example}

%%%%%%%%%%%%%%%%%%%%%%%%%%%%%%%%%%%%%%%%%%%%%%%%%%%%%%%%%%%%%%%%%%%%%%%%%%%%%%%%%%%%%%%%%%%%%%%%%%%%%%%%%%%%%%%%%%%%%%%%%%%%%%%%%%%%%%%%%%%%%%%%%%%%%%%%%%%%%%
\section{Bases and acyclic orientations}\label{section4}
In this section we describe a family of bases for the root spaces of a given Borcherds algebra under the assumptions of Section~\ref{section3}, namely that the coefficients of the real nodes are less or equal to one. 
\textit{Since $\lie g_{\eta(\bold k)}\neq 0$ implies that $\text{supp}(\bold k)=\{i\in I: k_i\neq 0\}$ is connected we will assume without loss of generality for the rest of this section that $I$ is connected and $I=\text{supp}(\bold k)$ and in particular $I$ is finite}. We freely use the notations introduced in the previous sections.

\subsection{}  
Let us first fix some notations. Let $I^*$ be the free monoid generated by $I$. Note that $I^*$ has a total order given by the lexicographical order. The free partially commutative monoid associated with $G$ is denoted by $M(I, G):=I^*/\sim$, where $\sim$ is generated by the relations $$ab\sim ba, \ \ (a, b)\notin E(G). $$
We associate to each element $a\in M(I,G)$ the unique element $\tilde{a}\in I^*$ which is the maximal element in the equivalence class of $a$ with respect to the lexicographical order. A total order on $M(I, G)$ is then given by
\begin{equation}\label{totor}a<b :\Leftrightarrow \tilde{a}<\tilde{b}.\end{equation}
Let $\bold w\in M(I, G)$ and write $\bold w=i_1\cdots i_r$. Define $|\bold w|=r$, $i(\bold w)=|\{j : i_j=i\}|$ for all $i\in I$ and $\text{supp}(\bold w)=\{i\in I :  i(\bold w)\neq 0\}$. The weight of $\bold w$ is denoted by 
$$\text{wt}(\bold w)=\sum\limits_{i\in I} i(\bold w)\alpha_i.$$
The initial alphabet of $\bold w$ is a multiset denoted by $\rm{IA_m}(\bold w)$ and defined by $i\in \rm{IA_m}(\bold w)$ (counted with multiplicities) if and only if $\exists \ \bold u\in M(I, G)$ such that $\bold w=\bold u i$. The underlying set is denoted by $\rm{IA}(\bold w)$. For example $\rm{IA_m}(1233)=\{3, 3\}$,  and $\rm{IA}(1233)=\{3\}$ for the complete graph $G=K_3$. The right normed Lie word associated with $\bold w$ is defined by 
$$e(\bold w)=[e_{i_1},[e_{i_{2}},[\cdots [e_{i_{r-1}}, e_{i_r}]]\in \lie g.$$
Using the Jacobi identity, it is easy to see that the association $\bold w\mapsto e(\bold w)$ is well defined.

\subsection{}The following lemma is straightforward.
\begin{lem}\label{strIA}
Let $\bold w=i_1\cdots i_r\in M(I,G)$. Then $|\rm{IA_m}(\bold w)|=1$ if and only if $\bold w$ satisfies the following condition:
$$\text{
given any $1\le k<r$  there  exists $k+1\le j\le r$ such that $(i_k, i_j)\in E(G)$}
$$
or equivalently $i_{r-1}\neq i_r$ and the subgraph generated by $i_k,\dots,i_r$ is connected for any $1\le k<r.$ 

\hfill\qed
\end{lem}
\subsection{} For $i\in I$, we introduce the so--called $i$--form of a given word $\bold w\in M(I, G)$.
\begin{prop}\label{iform}
Let $\bold w\in M(I, G)$ with $\rm{IA}(\bold w)=\{i\}$. Then there exists unique $\bold w_1,\dots,\bold w_{i(\bold w)}\in M(I, G)$ such that 
\begin{enumerate}
\item $\bold w=\bold w_1\cdots \bold w_{i(\bold w)}$
\item $\rm{IA_m}(\bold w_j)=\{i\}$,\ \ $i(\bold w_j)=1$ for each $1\le j\le i(\bold w)$.
\end{enumerate}
\end{prop}

\begin{pf}
We prove this result by induction on $i(\bold w).$ If $i(\bold w)=1$, there is nothing to prove; so assume that $i(\bold w)>1.$ We choose an expression $\bold w=i_1\cdots i_{r-1}i \in M(I, G)$  of $\bold w$ such that $i_k=i$, $i_{\ell}\neq i$ for all $k<\ell<r$ and $k$ is minimal with this property. To be more precise, if $\bold w=i_1'\cdots i_{r-1}'i$ is another expression of $\bold w$, with $i_{k'}=i$ and $i_{\ell}\neq i$ for all $k'<\ell<r$, then $k'\geq k$. We set $\bold w_{i(\bold w)}=i_{k+1}\cdots i_{r-1}i$. It is clear that $\rm{IA_m}(\bold w_{i(\bold w)})=\{i\}$ and $\bold w=\bold u \bold w_{i(\bold w)}$ where $\bold u=i_1\cdots i_{k}$. The minimality of $k$ implies that $\text{IA}(\bold u)=\{i\}$. Since $i(\bold u)=i(\bold w)-1$ we get by induction
$$\bold u=\bold w_1\cdots \bold w_{i(\bold w)-1}$$
such that $\rm{IA_m}(\bold w_j)=\{i\}$ for each $1\le j\le i(\bold w)-1$. 

Now we prove the uniqueness part; assume  that $\bold w=\bold w_1'\cdots \bold w_{i(\bold w)}'=\bold u' \bold w_{i(\bold w)}'$ is another expression such that $\rm{IA_m}(\bold w_j')=\{i\}$ for all $1\le j\le i(\bold w)$. Suppose $\bold w_{i(\bold w)}\neq \bold w_{i(\bold w)}'$ in $I^*$, which is only possible if there exists $i_p$ in $\bold u$, say $i_p\in \bold w_p$ with $p<i(\bold w)$ which we can pass through $i_{p+1},i_{p+2},\dots,i_{t-1},i_t=i\in \bold w_p$. This contradicts $|\rm{IA_m}(\bold w_p)|=1$. Hence $\bold w_{i(\bold w)}=\bold w_{i(\bold w)}'$ and the rest follows again by induction.
\end{pf}

The factorization of $\bold w$ in Proposition~\ref{iform} is called the $i$--form of $\bold w.$ \textit{In the rest of this section, we fix $i\in I$ and whenever we write $\bold w=\bold w_1\cdots \bold w_{i(\bold w)}$ we always assume that it is the $i$--form.} 
\subsection{}Here we will recall the combinatorics of Lyndon words and state the main theorem of this section; for more details about Lyndon words we refer the reader to \cite{Re93}. Consider the set $\mathcal{X}_{i}=\{\bold w\in M(I, G): \rm{IA}_m(\bold w)=\{i\}\}$ and recall that $\mathcal{X}_i$ (and hence $\mathcal{X}_i^*$) is totally ordered using \eqref{totor}. We denote by $FL(\mathcal{X}_i)$ the free Lie algebra generated by $\mathcal{X}_{i}$. A non--empty word $\bold w\in \mathcal{X}_i^{*}$ is called a Lyndon word if it satisfies one of the following equivalent definitions:
\begin{itemize}
\item $\bold w$ is strictly smaller than any of its proper cyclic rotations
\item $\bold w \in \mathcal{X}_i$ or $\bold w=\bold u \bold v$ for Lyndon words $\bold u$ and $\bold v$ with $\bold u<\bold v$.
\end{itemize}
There may be more than one choice of $\bold u$ and $\bold v$ with $\bold w=\bold u \bold v$ and $\bold u<\bold v$ but if $\bold v$ is of maximal possible length we call it the standard factorization.

To each Lyndon word $\bold w\in \mathcal{X}_{i}^*$ we associate a Lie word $L(\bold w)$ in $FL(\mathcal{X}_i)$ as follows. If $\bold w\in \mathcal{X}_i$, then $L(\bold w)=\bold w$ and otherwise $L(\bold w)=[L(\bold u),L(\bold v)]$, where $\bold w=\bold u \bold v$ is the standard factorization of $\bold w$.
The following result can be found in \cite{Re93} and is known as the Lyndon basis for free Lie algebras.
\begin{prop}\label{reut}
The set $\{L(\bold w): \bold w\in \mathcal{X}_i^* \text{ is a Lyndon word}\}$ forms a basis of $FL(\mathcal{X}_i)$.
\hfill\qed
\end{prop}

Let $\lie g^{i}$ the Lie subalgebra of $\lie g$ generated by $\{e(\bold w): \bold w\in \mathcal{X}_{i}\}$. By the universal property of $FL(\mathcal{X}_i)$ we have a surjective homomorphism
\begin{equation}\label{homsursp}\Phi: FL(\mathcal{X}_{i})\to \lie g^{i},\ \ \bold w\mapsto e(\bold w)\ \ \forall\ 
\bold w\in \mathcal{X}_{i}.\end{equation}
Using Proposition~\ref{reut} we immediately get that the image of $\Phi$ generates $\lie g^{i}$. It is natural to ask if in fact this procedure gives a basis for the root space $\lie g_{\eta(\bold k)}$; the main theorem of this section gives an answer to this question. Set 
$$C^{i}(\bold k, G)=\{\bold w\in \mathcal{X}_i^*:  \bold w \text{ is a Lyndon word},\  \rm{wt}(\bold w)=\eta(\bold k)\},\ \ \iota(\bold w)=\Phi\circ L(\bold w).$$
\begin{thm}\label{mainthmb}
The set $\left\{\iota(\bold w): \bold w\in C^{i}(\bold k, G)\right\}$ is a basis of the root space $\lie g_{\eta(\bold k)}$. Moreover, if $k_i=1$, the set
$$\left\{e(\bold w): \bold w\in \mathcal{X}_i,\ \rm{wt}(\bold w)=\eta(\bold k)\right\}$$
forms a right--normed basis of $\lie g_{\eta(\bold k)}$ and
$$C^{i}(\bold k, G)=\{\bold w\in \mathcal{X}_i:  e(\bold w)\neq 0,\  \rm{wt}(\bold w)=\eta(\bold k)\}.$$ 
\end{thm}
So if $k_i=1$, the above theorem implies that the root space $\lie g_{\eta(\bold k)}$ has a very special type of basis, namely that the right--normed Lie word $e(\bold w)$ is either zero or a basis element.

\subsection{}Here we proof Theorem~\ref{mainthmb}.
We set
$$\widetilde{B}^{i}(\bold k, G)=\left\{\bold w\in M(I, G) : \rm{wt}(\bold w)=\eta(\bold k) \ \rm{and} \ \ \rm{IA}(\bold w)=\{i\}\right\}.$$
Let $\bold w\in \widetilde{B}^{i}(\bold k, G)$ and write $\bold w=\bold w_1\cdots \bold w_{k_i}\in \mathcal{X}_i^{*}$. We say $\bold w$ is aperiodic if the elements in the cyclic rotation class of $\bold w$ are all distinct, i.e. all elements in 
$$C(\bold w)=\{\bold w_1\cdots \bold w_{k_i}, \bold w_2\cdots \bold w_{k_i}\bold w_1, \cdots, \bold w_{k_i}\bold w_1\cdots \bold w_{k_i-1}\}$$
are distinct. We naturally identify $C^{i}(\bold k, G)$ with the set
$$B^{i}(\bold k, G)=\{\bold w\in \widetilde{B}^{i}(\bold k, G): \bold w \text{ is aperiodic}\}/\sim,$$
where $\bold w\sim \bold w' \Leftrightarrow C(\bold w)=C(\bold w')$.
The following proposition is crucial for the proof of Theorem~\ref{mainthmb}. 
\begin{prop}\label{helppropp}
We have 
\begin{enumerate}
\item[(i)] The root space $\lie g_{\eta(\bold k)}$ is contained in $\lie g^i$.
\item[(ii)] Let $\bold w\in M(I,G)$ and $\rm{wt}(\bold w)=\eta(\bold k)$. Then 

$$e(\bold w)\neq 0\Longleftrightarrow \rm{IA}_m(\bold w)=\{i\}.$$
\item[(iii)]  We have
$$\text{mult } \eta(\bold k)= |B^{i}(\bold k, G)|.$$
\end{enumerate}
\end{prop}
The proof of the above proposition is postponed to the next subsection. We first show how this proposition proves Theorem~\ref{mainthmb}. Since $\lie g_{\eta(\bold k)}$ is contained in $\lie g^i$ we get with  Proposition~\ref{reut} and \eqref{homsursp} that $\left\{\iota(\bold w): \bold w\in C^{i}(\bold k, G)\right\}$ is a spanning set for $\lie g_{\eta(\bold k)}$ of cardinality equal to $|C^{i}(\bold k, G)|$. Therefore Proposition~\ref{helppropp} (iii) shows that this is in fact a basis. So in the special case when $k_i=1$ we get that $\left\{e(\bold w): \bold w\in \mathcal{X}_i,\ \rm{wt}(\bold w)=\eta(\bold k)\right\}$ is a basis. In order to finish the theorem we have to observe when a Lie word $e(\bold w)$ is non--zero, which is exactly answered by part (ii) of the proposition.
\subsection{Proof of Proposition~\ref{helppropp}(i)} 
This part of the proof is an easy consequence of the Jacobi identity.
\begin{lem}\label{spanlem}
Fix an index $i\in I$. Then the root space $\lie g_{\eta(\bold k)}$ is spanned by all right normed Lie words $e(\bold{w})$, where $\bold{w}\in M(I, G)$ is such that $\rm{wt}(\bold w)=\eta(\bold k)$ and $\rm{IA_m}(\bold w)=\{i\}.$
\begin{pf}
We fix $\bold{w}=i_1\cdots i_r$ and show that any element
$$e(\bold w,k)=\big[[[[e_{i_1},e_{i_2}],e_{i_3}]\cdots,e_{i_k}],[e_{i_{k+1}},[e_{i_{k+2}},\dots[e_{i_{r-1}}, e_{i_r}]]\big],\ 0\leq k <r$$
can be written as a linear combination of right normed Lie words $e(\bold{w}')$ with $\bold{w}'=j_1\cdots j_{r-1}i$. The claim finishes the proof since $e(\bold w)=e(\bold w,0)$ and $e(\bold w')\neq 0$ only if $\rm{IA_m}(\bold w')=\{i\}$. If $k=r-1$ this follows immediately using repeatedly $[x,y]=-[y,x]$ for all $x,y\in \lie g$. If $k<r-1$ we get from the Jacobi identity 
\begin{align}\label{asabov}e(\bold w,k)\notag &=e(\bold w,k+1)+\Big[e_{i_{k+1}},\big[[[[e_{i_1},e_{i_2}],e_{i_3}]\cdots,e_{i_k}],[e_{i_{k+2}},\dots[e_{i_{r-1}}, e_{i_r}]]\big]\Big]&\\&=e(\bold w,k+1)+\big[e_{i_{k+1}},e(\widetilde{\bold w},k+1)\big]
\end{align}
where $\widetilde{\bold w}=(i_1,\dots,\widehat{i_{k+1}},\cdots,i_r)$. An easy induction argument shows that each term in \eqref{asabov} has the desired property. 
\end{pf}
\end{lem}

\subsection{Proof of Proposition~\ref{helppropp}(ii)} Now we want to analyze when a right normed Lie word $e(\bold w)$ is non--zero. One has for $\bold w\in M(I, G)$, $e(\bold w)\neq 0$ implies $|\rm{IA_m}(\bold w)|=1$. Indeed we prove that the converse is also true.

\begin{lem}\label{keyprop}
The right normed Lie word $e(\bold{w})$ with $\rm{wt}(\bold w)=\eta(\bold k)$ is non--zero if and only if $|\rm{IA_m}(\bold w)|=1$.
\begin{pf}
Using Lemma~\ref{strIA} we prove that, for $\bold w=i_1\cdots i_r\in M(I, G)$,
the right normed Lie word $e(\bold w)$ is non--zero 
if and only if $\bold w$ satisfies the following condition:
\begin{equation}\label{neglem0}
\text{given any $1\le k<r$  there  exists
 $k+1\le j\le r$ such that $[e_{i_k}, e_{i_j}]\neq 0$.}
 \end{equation}

The remaining direction will be proven by induction on $r$, where the initial step $r=2$ obviously holds. So assume that $r>2$ and $e(\bold w)$ satisfies \eqref{neglem0}. 
We choose $p\in\{1,\dots,r-1\}$ be minimal such that $i_j\neq i_p$ for all $j>p$. Note that $p$ exists, since $i_{r-1}\neq i_r$. Let $I(p)=\{p_1,\dots,p_{k_{i_p}}\}$ be the elements satisfying $i_p=i_{p_j}$ for $1\leq j\leq k_{i_p}$. Note that $p_j\neq p$ implies $p_j<p$ by the choice of $p$. For a subset $S\subset I(r)$ let $\bold w(S)$ be the tuple obtained from $\bold w$ by removing the vertices $i_{q}$ for $q\in S$. 
The proof considers two cases.
\vspace{0,1cm}

\textit{Case 1:} We assume that $p<r-1$. We first show that $e(\bold w(S))$ satisfies \eqref{neglem0}. If \eqref{neglem0} is violated, 
we can find a vertex $i_k$ such that $k<p$, $i_k\neq i_p$ and $[e_{i_k},e_{i_\ell}]=0$ for all $\ell>k$ with $i_\ell\neq i_p$. 
We choose $k$ maximal with that property. By the minimality of $p$ we get $i_k\in\{i_{k+1},\dots, i_r\}\backslash\{i_p\}\supset \{i_{r-1},i_r\}$, 
say $i_k=i_m$ for some $k+1\leq m\leq r$. If $m>p$, we get $[e_{i_k},e_{i_t}]=[e_{i_m},e_{i_t}]\neq 0$ for some $k+1\leq t\leq r$ with $i_t\neq i_p$,
since $e(\bold w)$ satisfies \eqref{neglem0} and $p<r-1$. This is not possible by the choice of $i_k$. Hence $m<p$ and the maximality of $k$ implies the existence of a vertex $i_t$ such that $t>m$, $i_t\neq i_p$ and $[e_{i_k},e_{i_t}]=[e_{i_m},e_{i_t}]\neq 0$, which is once more a contradiction. Hence $e(\bold w(S))$ satisfies \eqref{neglem0} and is therefore non--zero by induction.

We consider the element $(\text{ad }f_{i_p})^{k_{i_p}}e(\bold w)$. If $i_p$ is a real node, we have $k_{i_p}=1$ and therefore 
$$(\text{ad }f_{i_p})e(\bold w)=-(\alpha_{i_{p}+1}+\cdots+\alpha_{i_r})(h_{i_p})e(\bold w(\{p\}))=(a_{i_p,i_{p}+1}+\cdots+a_{i_p,i_r})e(\bold w(\{p\})).$$
Since $(a_{i_p,i_{p}+1}+\cdots+a_{i_p,i_r})<0$ ($e(\bold w)$ satisfies \eqref{neglem0}) and $e(\bold w(\{p\}))\neq 0$ by the above observation we must have $e(\bold w)\neq 0$. If $i_p$ is not a real node then we have $a_{i_p,s}\leq 0$ for all $s\in I$. We get
$$(\text{ad }f_{i_r})^{k_{i_r}}e(\bold w)=C e(\bold w(I(r))),$$
for some non-zero constant $C$. Again we deduce $e(\bold w)\neq 0$.
\vspace{0,1cm}

\textit{Case 2:} We assume that $p=r-1$. In this case the rank of our Lie algebra is two. So $i_j\in I=\{1,2\}$ for all $1\leq j\leq r$. If $I^{\mathrm{re}}\neq \emptyset$, say $1\in I^{\mathrm{re}}$, we can finish the proof as follows. If $k_2=1$, there is nothing to prove. Otherwise, since $e(\bold w)$ satisfies \eqref{neglem0} we must have (up to a sign) $e(\bold w)=(\text{ad }e_{2})^{k_{2}}e_1$. Now $k_2>1$ implies $2\in I^{\mathrm{im}}$ and thus $\lie g_{\eta(\bold k)}\neq 0$ which forces $e(\bold w)\neq 0$. So it remains to consider the case when $I^{\mathrm{re}}=\emptyset$. In this case $\mathfrak{n}^+$ is the free Lie algebra generated by $e_1,e_2$ and the lemma is proven.
\end{pf}
\end{lem}

\subsection{}The rest of this section is dedicated to the proof of Proposition~\ref{helppropp}(iii). First we prove that $|\widetilde{B}^{i}(\bold k, G)|$ satisfies a recursion relation which is similar to the one in \eqref{recursss}. More precisely,

\begin{prop}\label{recursionforbi}
	We have
	$$|\widetilde{B}^{i}(\bold k, G)|=\sum\limits_{\ell|\bold k}\frac{k_i}{\ell}\left|B^{i}\left(\frac{\bold k}{\ell}, G\right)\right|.$$
\end{prop}
\begin{pf}
There exists a subset $\widehat{B}^{i}\left(\bold k, G\right)$ of $\widetilde{B}^{i}\left(\bold k, G\right)$ such that 
$$\{\bold w\in \widetilde{B}^{i}(\bold k, G): \bold w \text{ is aperiodic}\}=\dot{\bigcup_{\bold w\in \widehat{B}^{i}\left(\bold k, G\right)} }C(\bold w).$$
We clearly have $|\widehat{B}^{i}\left(\bold k, G\right)|=|B^{i}\left(\bold k, G\right)|$.
 Let $\ell\in \mathbb{Z}_{+}$ such that $\ell|\bold k$ and let $k=\frac{k_i}{\ell}.$ We consider the map
$\Phi_{\ell}: \widehat{B}^{i}\left(\frac{\bold k}{\ell}, G\right)\rightarrow \widetilde{B}^{i}(\bold k, G)$ defined by
 $$\bold w=\bold w_1\cdots \bold w_k\mapsto \underbrace{(\bold w_1\cdots \bold w_k)\cdots (\bold w_1\cdots \bold w_k)}_{\ell-\rm{times}}.$$ Since $\bold w$ is aperiodic it follows from the uniqueness of the $i$--form that $\Phi_{\ell}(\bold w)$ has exactly $k$ distinct elements in its cyclic rotation class. Choose another element $\bold w'\in \widehat{B}^{i}\left(\frac{\bold k}{\ell'}, G\right)$ such that $C(\Phi_{\ell}(\bold w))\cap C(\Phi_{\ell'}(\bold w'))\neq \emptyset$. We set $k'=\frac{k_i}{\ell'}$ and assume without loss of generality that $k\geq k'$; say $k=pk'+t$ with $0<t\leq k'$, $p\in \mathbb{Z}_+$. Let $\bold v= \bold v_1\cdots \bold v_k\in C(\bold w)$ and $\bold v'=\bold v_1'\cdots \bold v_{k'}'\in C(\bold w')$ such that $$\underbrace{\bold v\bold v\cdots \bold v}_{\ell-\rm{times}}=\underbrace{\bold v'\bold v'\cdots \bold v'}_{\ell'-\rm{times}}\in C(\Phi_{\ell}(\bold w))\cap C(\Phi_{\ell'}(\bold w')).$$ By the uniqueness property we get $\bold v_1\cdots \bold v_k=\bold v_{k-t+1}\cdots \bold v_k \bold v_1\cdots \bold v_{k-t}$. Since $\bold v$ is aperiodic we must have $t=k\leq k'$ and thus $\ell=\ell'$. Using once more the uniqueness property we obtain that $C(\bold w)=C(\bold w')$ and hence $\bold w=\bold w'$. It follows 
$$|\widetilde{B}^{i}(\bold k, G)|\geq \sum\limits_{\ell|\bold k}\frac{k_i}{\ell}\left|B^{i}\left(\frac{\bold k}{\ell}, G\right)\right|.$$ 
Now we prove that any element in $\bold w\in \widetilde{B}^{i}(\bold k, G)$ can be obtained by this procedure, i.e. we have to show that there exists $\ell\in\mathbb{Z}_+$ with $\ell | \bold k$ such that $\text{Im}(\Phi_{\ell})\cap C(\bold w)\neq \emptyset$. In what follows we construct an element in the aforementioned intersection. Let $\bold w=\bold w_1\cdots \bold w_{k_i}\in \widetilde{B}^{i}(\bold k, G)$ arbitrary. Note that the statement is clear if $\bold w$ is aperiodic. So let $1\leq r<s\leq k_i$ such that
$$\bold w_r\cdots \bold w_{k_i}\bold w_1\cdots \bold w_{r-1}=\bold w_s\cdots \bold w_{k_i}\bold w_1\cdots \bold w_{s-1}$$
and suppose that $k:=s-r$ is minimal with this property. Write $k_i=\ell k+t$, $0\leq t<k$. Again by the uniqueness of the $i$--form we obtain 
\begin{align*}\hspace{3cm} &\bold w_1=\bold w_{k+1}=\cdots=\bold w_{\ell k+1}=\bold w_{k_i-k+1},&\\&\bold w_2=\bold w_{k+2}=\cdots=\bold w_{\ell k+2}=\bold w_{k_i-k+2},&\\&\vdots &\\& \bold w_t=\bold w_{k+t}=\cdots=\bold w_{\ell k+t}=\bold w_{k_i-k+t},
&\\&\bold w_{t+1}=\bold w_{k+t+1}=\cdots=\bold w_{(\ell-1)k+t+1}=\bold w_{k_i-k+t+1},
&\\&\vdots&\\& \bold w_k=\bold w_{2k}=\cdots=\bold w_{\ell k}=\bold w_{k_i}.\end{align*}
Thus $\bold w_1=\cdots \bold w_{k_i}=\bold w_{t+1}=\cdots \bold w_{k_i}\bold w_1\cdots \bold w_{t}$ and the minimality of $k$ implies $t=0$. Therefore,
$$\bold w=\underbrace{(\bold w_1\cdots \bold w_k)\cdots (\bold w_1\cdots \bold w_k)}_{\ell-\rm{times}}.$$
If $\bold w_1\cdots \bold w_k$ is aperiodic we can choose the unique element in $\widehat{B}^{i}\left(\frac{\bold k}{\ell}, G\right)\cap C(\bold w_1\cdots \bold w_k)$ whose image is obviously an element of $\text{Im}(\Phi_{\ell})\cap C(\bold w)$. If not, we continue the above procedure with $\bold w_1\cdots \bold w_k$. This completes the proof.
\end{pf}
\begin{rem}
From the proof of Proposition~\ref{recursionforbi} we get that any element of $B^{i}(\bold k, G)$ is aperiodic if $k_i,i\in I$ are relatively prime.
\end{rem}
\subsection{} 
In this subsection we prove Theorem~\ref{mainthmb} for the special case $\bold k=\bold 1=(k_i=1 : i\in I)$ which will be needed later to prove Proposition~\ref{helppropp}(iii). We need the notion of an acyclic orientation for this. An acyclic orientation of a graph is an assignment of a direction to each edge that does not form any directed cycle. A node $j$ is said to be a sink of an orientation if no arrow is incident to $j$ at its tail. 
We denote by $O_{i}(G)$ the set of acyclic orientations of $G$ with unique sink $i$. Greene and Zaslavsky \cite[Theorem 7.3]{GZ83} proved a connection between the number of acyclic orientations with unique sink and the chromatic polynomial.  In particular, up to a sign
\begin{equation}\label{acyclicor}
|O_{i}(G)|= \text{ the linear coefficient of the chromatic polynomial of $G$}.
\end{equation}
Now we are able to prove,
\begin{prop}\label{k=1case}
The set 
$$\{e(\bold w) : \bold w\in \widetilde{B}^{i}\left(\bold 1, G\right)\}$$
is a basis of the root space $\lie g_{\eta(\bold 1)}$.
\end{prop}
\begin{pf}
 Given any $\bold w\in \widetilde{B}^{i}(\bold 1, G)$, we associate an acyclic orientation $G^{\bold w}$ as follows. Draw an arrow from the node $u$ to $v$ ($u\neq v$) if and only if $u$ appears on the left of $v$ in $\bold w$. 
 Note that the assignment $\bold w\mapsto G^{\bold w}$ is well defined, i.e. $G^{\bold w}$  does not depend on the choice of an expression of $\bold w$. It is easy to see that the node $i$ is the unique sink for the acyclic orientation $G^{\bold w}$. Thus the assignment $\bold w\mapsto G^{\bold w}$ defines a well--defined injective map $\widetilde{B}^{i}(\bold 1, G)\rightarrow O_{i}(G)$. Since $\{e(\bold w): \bold w\in \widetilde{B}^{i}(\bold 1, G)\}$ spans $\lie g_{\eta(\bold 1)}$ we obtain together with Corollary~\ref{recursionmult}
$$|O_{i}(G)|\geq |\widetilde{B}^{i}(\bold 1, G)|\geq \text{dim } \lie g_{\eta(\bold 1)}=|\pi^G_\bold 1(q)[q]|.$$
Now the proposition follows from \eqref{acyclicor}.
\end{pf}
As a corollary we get,
\begin{cor}
The number of acyclic orientations of $G$ with unique sink $i\in I$ is equal to $\text{mult } \eta(\bold 1)$.
\hfill\qed
\end{cor}
\subsection{}\label{clique}
We construct a new graph associated to the pair $(G, \bold k)$ which will help to determine the cardinality of $\widetilde{B}^i(\bold k,G)$. For each $j\in I$, take a clique (complete graph) of size $k_j$ with vertex set $\{j^1,\dots, j^{k_j}\}$ and join all vertices of the $r$--th and $s$--th clique if $(r,s)\in E(G).$ The resulting graph is called the join of $G$ with respect to $\bold k$, denoted by $G(\bold k)$. 

We consider  
$$M(i_1,\dots,i_{k_i}):=\widetilde{B}^{i^1}(\bold 1, G(\bold k))\ \dot{ \cup }\cdots \dot{ \cup }\ \widetilde{B}^{i^{k_i}}(\bold 1, G(\bold k))$$
and define a map $\phi:M(i_1,\dots,i_{k_i})\rightarrow \widetilde{B}^{i}(\bold k, G)$ by sending the vertices of the $j$--th clique $\{j^1,\dots, j^{k_j} \}$ to the vertex $j$ for all $j\in I$. This map is clearly surjective.
Let $\bold w\in M(i_1,\dots,i_{k_i})$ and let $\sigma \in S_{k_j}$, $j\in I$ which acts on $\bold w$ by permuting the entries in the $j$--th clique. It is easy to show that $|\rm{IA_m}(\sigma(\bold w))|=1$ and thus $\sigma(\bold w)\in M(i_1,\dots,i_{k_i})$. Therefore there is a natural action of the symmetric group $$\bold{S}_{\bold k}=\prod\limits_{j\in I} S_{k_j} \ \ \  \text {on}  \ \ \
M(i_1,\dots,i_{k_i}).$$ 
It is clear that this action induces a bijective map, 
$$\overline{\phi}: M(i_1,\dots,i_{k_i})/\bold S_{\bold k}\to \widetilde{B}^{i}(\bold k, G).$$
We will conclude our discussion by proving that the above action is free. Let $\sigma\in \bold S_{\bold k}$ and $\bold w\in M(i_1,\dots,i_{k_i})$ such that $\sigma \bold w=\bold w$. If $\sigma$ is not the identity we can find $j\in I$ and $1\leq \ell < s \leq k_j$ such that $\sigma(\ell)>\sigma(s)$. So if $j^\ell$ appears on the left of $j^{k}$ in $\bold w$, then $j^k$ appears on the left of $j^{\ell}$ in $\sigma(\bold w)$. Since $\sigma(\bold w)=\bold w$, this is only possible if $(j^{\ell},j^{k})\in E(G(\bold k))$ which is a contradiction.
Thus we have proved that 
\begin{equation}\label{einew}|\widetilde{B}^{i}(\bold k, G)|=\sum_{r=1}^{k_i}\frac{1}{\bold k!}|\widetilde{B}^{i^r}(\bold 1, G(\bold k))|.\end{equation}
\subsection{} We complete the proof of Proposition~\ref{helppropp}(iii). We have
\begin{align*}
&& \sum\limits_{\ell|\bold k}\frac{1}{\ell} |B^i\left(\bold k, G\right)|&=\frac{1}{k_i}|\widetilde{B}^{i}(\bold k, G)|&& \text{by Proposition~\ref{recursionforbi}}\\
&& &=\frac{1}{k_i}\sum_{r=1}^{k_i}\frac{1}{\bold k!} |\widetilde{B}^{i^r}(\bold 1, G(\bold k))|&& \text{by \eqref{einew}}\\
&& &=\frac{|\pi_{\bold 1}^{G(\bold k)}(q)[q]|}{\bold k!}&& \text{by Proposition~\ref{k=1case}}\\
&& &=|\pi^{G}_{\bold k}(q)[q]|&& \text{}\\
&& &=\sum\limits_{\ell|\bold k}\frac{1}{\ell}\text{mult }\eta\left(\frac{\bold k}{\ell}\right)&& \text{by \eqref{recursss}.}
\end{align*}
Now induction on g.c.d of $\bold k$ implies that $|B^i\left(\bold k, G\right)|=\text{mult } \eta(\bold k)$ for all $\bold k$.
\begin{example}We finish this section by determining a basis for the root space $\lie g_{\eta(\bold k)}$ where $G$ is the graph defined in Example~\ref{multicoloring} and $\mathbf{k}=(2,1,1,1)$. The following is a list of acyclic orientations of $G(\bold k)$ with unique sink $i=2$

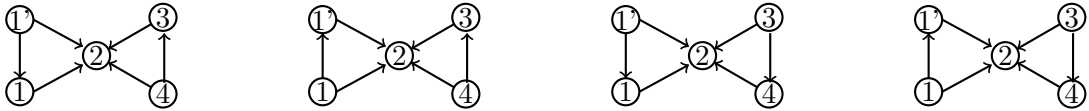
\vspace{0,5cm}

\hspace{-4,5cm}
\begin{tikzpicture}[scale=.6]
    \draw (-1,0) node[anchor=east]{};
		\draw[xshift=3 cm,thick] (33 mm, 8 mm) circle(.3cm) node[]{1'};
    \draw[xshift=3 cm,thick] (33 mm, -8 mm) circle(.3cm) node[]{1};
    \draw[xshift=4 cm,thick] (4 cm,0) circle(.3cm)node[]{2};
    \draw[xshift=8 cm,thick] (30: 17 mm) circle (.3cm)node[]{3};
    \draw[xshift=8 cm,thick] (-30: 17 mm) circle (.3cm)node[]{4};
  %  \draw[dotted,thick] (0.3 cm,0) -- +(1.4 cm,0);
    \foreach \y in {3.15,...,3.15}
		\draw[xshift=3 cm,thick,->] (33 mm, 5 mm) -- (33 mm, -5 mm);
    \draw[xshift=3 cm,thick,->] (36 mm, 8 mm) -- (47 mm,2 mm);
		\draw[xshift=3 cm,thick,->] (36 mm, -8 mm) -- (47 mm,-2 mm);
    \draw[xshift=8 cm,thick, <-] (30: 3 mm) -- (30: 14 mm);
    \draw[xshift=8 cm,thick,<-] (-30: 3 mm) -- (-30: 14 mm);
      \draw[thick,->] (95mm, -6mm) -- +(0, 11mm);
\end{tikzpicture}
\hspace{-3cm}
\begin{tikzpicture}[scale=.6]
    \draw (-1,0) node[anchor=east]{};
		\draw[xshift=3 cm,thick] (33 mm, 8 mm) circle(.3cm) node[]{1'};
    \draw[xshift=3 cm,thick] (33 mm, -8 mm) circle(.3cm) node[]{1};
    \draw[xshift=4 cm,thick] (4 cm,0) circle(.3cm)node[]{2};
    \draw[xshift=8 cm,thick] (30: 17 mm) circle (.3cm)node[]{3};
    \draw[xshift=8 cm,thick] (-30: 17 mm) circle (.3cm)node[]{4};
  %  \draw[dotted,thick] (0.3 cm,0) -- +(1.4 cm,0);
    \foreach \y in {3.15,...,3.15}
		\draw[xshift=3 cm,thick,<-] (33 mm, 5 mm) -- (33 mm, -5 mm);
    \draw[xshift=3 cm,thick,->] (36 mm, 8 mm) -- (47 mm,2 mm);
		\draw[xshift=3 cm,thick,->] (36 mm, -8 mm) -- (47 mm,-2 mm);
    \draw[xshift=8 cm,thick, <-] (30: 3 mm) -- (30: 14 mm);
    \draw[xshift=8 cm,thick,<-] (-30: 3 mm) -- (-30: 14 mm);
      \draw[thick,->] (95mm, -6mm) -- +(0, 11mm);
\end{tikzpicture}
\hspace{-3cm}
\begin{tikzpicture}[scale=.6]
    \draw (-1,0) node[anchor=east]{};
		\draw[xshift=3 cm,thick] (33 mm, 8 mm) circle(.3cm) node[]{1'};
    \draw[xshift=3 cm,thick] (33 mm, -8 mm) circle(.3cm) node[]{1};
    \draw[xshift=4 cm,thick] (4 cm,0) circle(.3cm)node[]{2};
    \draw[xshift=8 cm,thick] (30: 17 mm) circle (.3cm)node[]{3};
    \draw[xshift=8 cm,thick] (-30: 17 mm) circle (.3cm)node[]{4};
  %  \draw[dotted,thick] (0.3 cm,0) -- +(1.4 cm,0);
    \foreach \y in {3.15,...,3.15}
		\draw[xshift=3 cm,thick,->] (33 mm, 5 mm) -- (33 mm, -5 mm);
    \draw[xshift=3 cm,thick,->] (36 mm, 8 mm) -- (47 mm,2 mm);
		\draw[xshift=3 cm,thick,->] (36 mm, -8 mm) -- (47 mm,-2 mm);
    \draw[xshift=8 cm,thick, <-] (30: 3 mm) -- (30: 14 mm);
    \draw[xshift=8 cm,thick,<-] (-30: 3 mm) -- (-30: 14 mm);
      \draw[thick,<-] (95mm, -6mm) -- +(0, 11mm);
\end{tikzpicture}
\hspace{-3cm}
\begin{tikzpicture}[scale=.6]
    \draw (-1,0) node[anchor=east]{};
		\draw[xshift=3 cm,thick] (33 mm, 8 mm) circle(.3cm) node[]{1'};
    \draw[xshift=3 cm,thick] (33 mm, -8 mm) circle(.3cm) node[]{1};
    \draw[xshift=4 cm,thick] (4 cm,0) circle(.3cm)node[]{2};
    \draw[xshift=8 cm,thick] (30: 17 mm) circle (.3cm)node[]{3};
    \draw[xshift=8 cm,thick] (-30: 17 mm) circle (.3cm)node[]{4};
  %  \draw[dotted,thick] (0.3 cm,0) -- +(1.4 cm,0);
    \foreach \y in {3.15,...,3.15}
		\draw[xshift=3 cm,thick,<-] (33 mm, 5 mm) -- (33 mm, -5 mm);
    \draw[xshift=3 cm,thick,->] (36 mm, 8 mm) -- (47 mm,2 mm);
		\draw[xshift=3 cm,thick,->] (36 mm, -8 mm) -- (47 mm,-2 mm);
    \draw[xshift=8 cm,thick, <-] (30: 3 mm) -- (30: 14 mm);
    \draw[xshift=8 cm,thick,<-] (-30: 3 mm) -- (-30: 14 mm);
      \draw[thick,<-] (95mm, -6mm) -- +(0, 11mm);
\end{tikzpicture}

So we have $\text{mult }\eta(\bold k)=2$ and if $i=2$, the following right--normed Lie words form a basis of $\lie g_{\eta(\bold k)}$
$$[e_4,[e_3,[e_1,[e_1,e_2]]]],\ [e_3,[e_4,[e_1,[e_1,e_2]]]].$$
Similarly if $i=1$ we have the basis
$$[e_1,[e_3,[e_4,[e_2,e_1]]],\ [e_1,[e_4,[e_3,[e_2,e_1]]].$$
\end{example}
%%%%%%%%%%%%%%%%%%%%%%%%%%%%%%%%%%%%%%%%%%%%%%%%%%%%%%%%%%%%%%%%%%%%
%%%%%%%%%%%%%%%%%%%%%%%%%%%%%%%%%%%%%%%%%%%%%%%%%%%%%%%%%%%%%%%%%%%%
%%%%%%%%%%%%%%%%%%%%%%%%%%%%%%%%%%%%%%%%%%%%%%%%%%%%%%%%%%%%%%%%%%%

\section{Hilbert series for partially commutative Lie algebras}\label{section5} 
In this section we will compute the Hilbert series of the $q$--fold tensor product of the universal enveloping algebra of free partially commutative Lie algebras.
{\em In what follows, we assume that the Borcherds algebra $\lie g$ has only imaginary simple roots.} In this case, it is easy to see that $\lie n^+$ is isomorphic to the free partially commutative Lie algebra associated to $G$.
We will use the ideas of the previous sections to show that the coefficients of the Hilbert series are determined by the generalized chromatic polynomials evaluated at negative integers. As an application, we will give a Lie theoretic proof of Stanley's reciprocity theorem of chromatic polynomials \cite{S73}. 
\subsection{} We first recall the definition of Hilbert series. 
Let $\Gamma$ be a semigroup with at most countably infinite elements and $\lie a$ be a $\Gamma$--graded Lie algebra with finite dimensional homogeneous spaces, i.e.
 $$\lie a=\bigoplus\limits_{\alpha\in \Gamma}\ \lie a_\alpha \ \ \ \ \text{and \ \ \  $\rm{dim}(\lie a_{\alpha})<\infty$ for all $\alpha\in \Gamma$}.$$ 
 The Hilbert series of $\lie a$ is defined as $H_{\Gamma}({\lie a})=\sum\limits_{\alpha\in \Gamma}^{}\text{dim } \lie a_{\alpha}\ e^{\alpha}.$ The $\Gamma$--grading of $\lie a$ induces a $\Gamma\cup \{0\}$--grading on the universal enveloping algebra $\mathbf{U}(\lie a)$ and similarly we define
 $$H_{\Gamma}(U(\lie a))=1+\sum\limits_{\alpha\in \Gamma} (\text{dim } \mathbf{U}(\lie a)_\alpha)\ e^{\alpha}.$$ 
The proof of the following proposition is standard using the Poincar{\'e}--Birkhoff--Witt theorem. 
 \begin{prop}\label{hilfsproph}
 The Hilbert series of $\mathbf{U}(\lie a)$ is given by
 $$H_{\Gamma}(\mathbf{U}(\lie a))= \frac{1}{\prod\limits_{\alpha \in \Gamma}\left(1 - e^{\alpha}\right)^{\dim \mathfrak{a}_{\alpha}}}.$$
\end{prop}
\subsection{}
Set $\Gamma=Q_+\backslash\{0\}.$ Since $\lie n^+$ is $\Gamma$--graded, the following statement is a straightforward application of the denominator identity:
\begin{cor}\label{corrrr}
The Hilbert series of $\mathbf{U}(\lie n^+)$ is given by
$$H_{\Gamma}(\mathbf{U}(\lie n^+)) =  \left(\sum_{ \gamma \in
	\Omega} (-1)^{\mathrm{ht}(\gamma)} e^{\gamma}\right)^{-1}.$$
\end{cor}
\begin{pf}
 Since all the simple roots of $\lie g$ are imaginary, the denominator identity \ref{denominator} (applied to $\lie n^+$) becomes 
 $$\sum_{\gamma \in
\Omega} (-1)^{\mathrm{ht}(\gamma)} e^{\gamma}  = \prod_{\alpha \in \Delta_+} (1 - e^{\alpha})^{\dim \mathfrak{g}_{\alpha}} 
$$ 
 Using Proposition~\ref{hilfsproph} we see that
 $$H_{\Gamma}(\mathbf{U}(\lie n^+)) =
 \frac{1}{\prod\limits_{\alpha \in \Gamma}\left(1 - e^{\alpha}\right)^{\dim \mathfrak{g}_{\alpha}}}=\frac{1}{\prod\limits_{\alpha \in \Delta_+}\left(1 - e^{\alpha}\right)^{\dim \mathfrak{g}_{\alpha}}}$$
where the last equality follows from $\Delta_+=\{\alpha\in Q_+\backslash\{0\}: \text{dim } \lie g_\alpha \neq 0 \}.$
 Now the result is immediate from the denominator identity.
\end{pf}
\subsection{}In this subsection, we show that the Hilbert series of the $q$--fold tensor product of the universal enveloping algebra $\mathbf{U}(\lie n^+)$ is determined by the evaluation of the generalized chromatic polynomial at $-q$. 
We naturally identify $\alpha\in Q_+$ with a tuple on non--negative integers. We prove, 
\begin{thm}\label{tensoruniversal}  Let $q\in \mathbb{N}$. Then the Hilbert series of $\mathbf{U}(\lie n^+)^{\otimes q}$ is given by
 $$H_{\Gamma}(\mathbf{U}(\lie n^+)^{\otimes q}) = \sum\limits_{\alpha\in Q_+}(-1)^{\rm{ht}(\alpha)} \pi_\alpha^G(-q)\ e^{\alpha}.$$
In particular,
$$\text{dim }(\mathbf{U}(\lie n^+)^{\otimes q})_{\alpha} = (-1)^{\rm{ht}(\alpha)} \pi_\alpha^G(-q),\ \ \text{ for all } \alpha\in Q_+.$$
\end{thm}
\begin{pf}
We obviously have $H_{\Gamma}(\mathbf{U}(\lie n^+)^{\otimes q})=H_{\Gamma}(\mathbf{U}(\lie n^+))^q$. From Proposition \ref{helprop}, we get
\begin{eqnarray}  \left(\sum_{ \gamma \in
	\Omega} (-1)^{\mathrm{ht}(\gamma)} e^{\gamma}\right)^{-q} & =& \sum\limits_{\alpha\in Q_+}(-1)^{\text{ht}(\alpha)} \pi_\alpha^G(-q)\ e^{\alpha}.
 \end{eqnarray}
 A straightforward application of Corollary~\ref{corrrr} finishes the proof.
 \end{pf}

\subsection{}The following is an immediate consequence of Theorem~\ref{tensoruniversal}. For a multiset $S$ we let $\alpha(S)=\sum_{i\in S}\alpha_i$.
 \begin{cor}\label{valueatq}
  Let $q\in \mathbb{N}$ and $\alpha=\sum_{i\in I}k_i\alpha_i\in Q^+$. Then the generalized chromatic polynomial of $G$ satisfies,
  $$ \pi^G_{\alpha}(-q)=\sum\limits_{(S_1,\dots, S_q)} \pi^G_{\alpha(S_1)}(-1)\cdots \pi^G_{\alpha(S_q)}(-1)$$
	where the sum runs over all ordered partitions $(S_1,\dots, S_q)$ of the multiset 
\begin{equation}\label{multiss}\{\underbrace{\alpha_i,\dots, \alpha_i}_{k_i\, \mathrm{times}}: i\in I\}.\end{equation}
 \end{cor}
\begin{pf}
 By Theorem~\ref{tensoruniversal}, we get $\text{dim }(\mathbf{U}(\lie n^+)^{\otimes q})_{\alpha} = (-1)^{\text{ht}(\alpha)} \pi_{\alpha}^G(-q)$. But
 $$(\mathbf{U}(\lie n^+)^{\otimes q})_{\alpha}= \bigoplus\limits_{(S_1,\dots, S_q)}\mathbf{U}(\lie n^+)_{\alpha(S_1)}\otimes \cdots \otimes \mathbf{U}(\lie n^+)_{\alpha(S_q)}$$
 where the sum runs over all ordered partitions $(S_1,\dots, S_q)$ of the multiset \eqref{multiss}. Again by Theorem~\ref{tensoruniversal} we know that $\text{dim }(\mathbf{U}(\lie n^+))_{\alpha(S_i)} = (-1)^{\text{ht}(\alpha(S_i))} \pi_{\alpha(S_i)}^G(-1)$ for all $S_i$. Putting all this together we get the desired result.
\end{pf}

\subsection{}We give another application of Theorem~\ref{tensoruniversal}, namely we show how this can be used to give a different proof of Stanley's reciprocity theorem of chromatic polynomials \cite[Theorem 1.2]{S73}. Note that the universal enveloping algebra of $\lie n^+$ is isomorphic 
to the free associative algebra generated by $e_i$, $i\in I$ with the relations $e_ie_j=e_je_i$ for all $(i, j)\notin E(G).$   
{\em Assume for the rest of this section that $G$ is a finite graph with index set $I$.} Suppose we are given an acyclic orientation $\mathcal{O}$ of $G$. Let
$$e_{\mathcal{O}}=e_{i_1}\cdots e_{i_n}$$
be the unique monomial in $\mathbf{U}(\lie n^+)$ determined by the following two conditions: 
$(i_1,\dots, i_n)$ \text{is a permutation of $I$} and if we have an arrow $i_r\rightarrow i_s$ in $\mathcal{O}$, then we require $r<s$.
Clearly, this gives a bijective correspondence between the set of acyclic orientations of $G$ and a vector space basis
of $\mathbf{U}(\lie n^+)_{\alpha(I)}$. 
For a map $\sigma: S\to \{1, 2,\dots,q\}$, we say $(\sigma, \mathcal{O})$ is a $q$--compatible pair if for each directed edge $i\rightarrow j$ in $\mathcal{O}$
we have $\sigma(i)\ge \sigma(j)$. 
\begin{thm}
 The number of $q$--compatible pairs of $G$ is equal to $(-1)^{|I|}\pi_{\alpha(I)}(-q)$. In particular, $(-1)^{|I|}\pi_{\alpha(I)}(-1)$ counts the number of acyclic orientations of $G$.
\end{thm}
\begin{pf}
 Let $(\sigma, \mathcal{O})$ be a $q$--compatible pair. For $\mathcal{O}$ we consider the unique monomial 
$$e_{\mathcal{O}}=e_{i_1}\cdots e_{i_n}.$$
 Set $S_{q-{i}+1}=\{v\in S: \sigma(v)=i\}$ for $1\le i\le q$ and let $\mathcal{O}_i$ the acyclic orientation of the subgraph spanned by $S_i$ which is induced by $\mathcal{O}$. Then $(S_1,\dots, S_q)$ is a partition of $S$ and since $(\sigma, \mathcal{O})$ is $q$--compatible pair, it is easy to see that
 $$e_{\mathcal{O}}=e_{\mathcal{O}_1}\cdots e_{\mathcal{O}_q}.$$
 
 On the other hand, suppose $(S_1,\dots, S_q)$ is a partition of $S$ and let $\mathcal{O}_i$ an acyclic orientation of the subgraph spanned by $S_i$ for $1\le i\le q$. We extend the acyclic orientations to an acyclic orientation of $G$. If $(i,j)\in E(G)$ such that $i\in S_r$ and $j\in S_k$ and $r<k$, then we give the orientation $i\rightarrow j$. This defines an acyclic orientation $\mathcal{O}$ of $G$. Further, define $\sigma: S\to \{1,\cdots, q\}$ by $\sigma(S_i)=q-i+1$ for all $1\le i\le q.$ Then it is easy to see that $(\sigma, \mathcal{O})$ is a $q$--compatible pair.
Thus we proved,   
 $$ \text{\# of $q$--compatible pairs of $G$}=\sum\limits_{(S_1,\dots, S_q)} (-1)^{|S_1|} \pi^G_{\alpha(S_1)}(-1)\cdots (-1)^{|S_q|}\pi^G_{\alpha(S_q)}(-1)$$where the sum runs over all $(S_1,\dots, S_q)$ partitions of $S.$
Now Corollary~\ref{valueatq} completes the proof.

\end{pf}

\subsection{} We consider the hight grading on $\lie n^+$ 
$$\lie n^+ = \bigoplus\limits_{k=1}^{\infty}\lie g_k,\ \ \ \  \ \lie g_k=\bigoplus\limits_{\overset{\text{ht } \alpha = k}{\alpha \in \Delta_{+}} }\lie g_\alpha$$
and show that in this case the Lucas polynomials show up as coefficients of the Hilbert series when the complement graph is triangle free. The hight grading naturally arises in the study of lower central series of right angled Artin groups. Let $\mathcal{G}$ be the right angled Artin group associated to $G$. The lower central series of $\mathcal{G}$ is the sequence of normal subgroups defined by 
$\mathcal{G}_1=\mathcal{G}$ and $\mathcal{G}_{n+1}=[\mathcal{G}_{n}, \mathcal{G}]$ for $n\ge 1.$ Then the associated graded Lie algebra of $\mathcal{G}$,
$$\rm{gr}_{\mathbb{C}}(\mathcal{G})= \bigoplus\limits_{k=1}^{\infty} \mathbb{C}\otimes_\mathbb{Z}(\mathcal{G}_{k+1}/\mathcal{G}_k)$$ is isomorphic to $\lie n^+$, see \cite[Theorem 3.4]{PS06}.
The Lie bracket $[x, y]$ of $\rm{gr}_{\mathbb{C}}(\mathcal{G})$  is induced from
the group commutator and the grading is given by bracket length. In particular, the abelian group $\mathcal{G}_{k+1}/\mathcal{G}_k$ is torsion free and its rank is equal to the dimensions of $\lie g_k$.

For the rest of this subsection we relate $M_k:=\text{dim } \lie g_k$ to the independent set polynomial of $G$ defined by 
$$I(X, G)=\sum\limits_{i\ge 0}c_i X^i $$
where $c_0=1$ and $c_i, i\geq 1$ is the number of independent subsets of G of size $i$. Using
$$U=\prod\limits_{k \ge 1} (1 - X^{k})^{M_k}=\sum_{j \ge 0} ~ (-1)^{j} ~ c_{j} ~ X^{j}$$
we get by comparing coefficients of $X^{k}$ on both sides in $-\text{log} ~ (U)$ that 
$$
N_{k} = \frac{1}{k} ~ \sum\limits_{d|k} ~ \frac{k}{d} ~ M_{\frac{k}{d}}= \frac{1}{k} ~ \sum_{d|k} ~ d ~ M_{d},
$$
where
\begin{eqnarray}
N_{k}:=(-1)^{k+1} c_{k} + \frac{1}{2} (-1)^{k+2} \left\{ \sum\limits_{\overset{j_{1} + j_{2} =k} {j_{1},j_{2} \ge 1}} c_{j_{1}} c_{j_{2}} \right\} + 
\frac{1}{3} (-1)^{k+3} \left\{ \sum\limits_{\overset{j_{1} + j_{2} + j_{3} = k} {j_{1},j_{2}, j_{3} \ge 1}} c_{j_{1}} c_{j_{2}} c_{j_{3}} \right\} + \cdots 
\end{eqnarray}
Using the M\"{o}bius inversion formula we recover a result of \cite[Theorem III.3.]{DK92}, 
\begin{equation}\label{ranksbas} M_{k} =  \sum\limits_{d|k} ~ \frac{\mu(d)}{d} ~ N_{\frac{k}{d}}.\end{equation}

We remark that the connection between the numbers $N_k$'s and the coefficients of the independent set polynomial of $G$ is not explained in \cite{DK92}.

\subsection{} In this subsection, we consider a graph $G$ such that  $G^{\rm{c}}$ (the complement graph of $G$) is triangle free and determine explicitly the numbers $N_k's$. In this case, these numbers are related to the Lucas polynomials $\big<\ell\big>_{s, t}$ in two variables $s$ and $t$. The Lucas polynomial is defined by the recurrence relation $$\big<\ell\big>_{s,t}=s\big<\ell-1\big>_{s,t}+t\big<\ell-2\big>_{s,t} \ \ \text{for} \ \ \ell\ge 2$$ with the initial values $\big<0\big>_{s,t}=2$ and $\big<1\big>_{s,t}=s.$ 
There is a closed formula (see \cite[Proposition 2.1]{ACMS14})
$$\big<\ell\big>_{s,t}=\sum\limits_{j\ge 0}\frac{\ell}{\ell-j}{\ell-j \choose j }
t^j s^{\ell-2j}.$$
We show,
\begin{prop} Let $G$ be a graph such that $G^{\rm{c}}$ is triangle free. Let $v$ be the number of vertices of $G^{\rm{c}}$ and $e$ the number of edges of $G^{\rm{c}}.$ Then for $k\ge 1$, we have

$$N_k=\frac{1}{k}\big<k\big>_{v, -e},\ \ M_k = \frac{1}{k}\sum\limits_{d|k} ~ \mu\left(\frac{k}{d}\right) \big<d\big>_{v, -e}.$$
\proof
We have $I(-X, G)=1-vX+eX^2$. Hence applying $-\text{log}$ we get,
$$\sum\limits_{k=1}^{\infty}\frac{(vX)^k}{k}(1-(e/v)X)^k=\sum\limits_{k=1}^{\infty}\frac{1}{k}
\left(\sum\limits_{r=0}^{k}(-1)^r {k \choose r}e^rv^{k-r}X^{k+r} \right).$$ 
The coefficient of $X^{k}$ in $-\text{log}~(1-vX+eX^2)$ is 
$$N_k=\sum\limits_{\overset{p+r=k}{p\ge r\ge 0}}\frac{(-1)^r}{p} {p \choose r} e^r v^{p-r}=\sum\limits_{j\ge 0}\frac{(-1)^j}{k-j}{k-j \choose j }
e^j v^{k-2j}=\frac{1}{k}\big<k\big>_{v, -e}.$$
The rest follows from \eqref{ranksbas}.

\endproof
\end{prop}

%%%%%%%%%%%%%%%%%%%%%%%%%%%%%%%%%%%%%%%%%%%%%%%%%%%%%%%%%%%%%%%%%%%%%%%%%%%%%%%%%%%%%%%%%%%%%%%%%%%%%%%%%%%%%%%%%%%%%%%%%%%%%%%%%%%%%%%%

%%%%%%%%%%%%%%%%%%%%%%%%%%%%%%%%%%%%%%%%%%%%%%%%%%%%%%%%%%%%%%%%%%%%%%%%%%%%%%%%%%%%%%%%%%%%%%%%%%%%%%%%%%%%%%%%%%%%%%%%%%%%%%%%%%%%%%%%

%%%%%%%%%%%%%%%%%%%%%%%%%%%%%%%%%%%%%%%%%%%%%%%%%%%%%%%%%%%%%%%%%%%%%%%%%%%%%%%%%%%%%%%%%%%%%%%%%%%%%%%%%%%%%%%%%%%%%%%%%%
%%%%%%%%%%%%%%%%%%%%%%%%%%%%%%%%%%%%%%%%%%%%%%%%%%%%%%%%%%%%%%%%%%%%%%%%%%%%%%%%%%%%%%%%%%%%%%%%%%%%%%%%%%%%%%%%%%%%%%%%%%

%%%%%%%%%%%%%%%%%%%%%%%%%%%%%%%%%%%%%%%%%%%%%%%%%%%%%%%%%%%%%%%%%%%%%%%%%%%%%%%%%%%%%%%%%%%%%%%%%%%%%%%%%%%%%%%%%%%%%%%%%
% References
%%%%%%%%%%%%%%%%%%%%%%%%%%%%%%%%%%%%%%%%%%%%%%%%%%%%%%%%%%%%%%%%%%%
\bibliographystyle{plain}
\bibliography{kv-bib}

\end{document}